\documentclass[reqno,final,a4paper]{amsart}
\usepackage[utf8]{inputenc}
\usepackage[notref,notcite]{showkeys} 
\usepackage[shortlabels]{enumitem}
\setitemize{itemsep=0.3em}
\setenumerate{itemsep=0.3em}
\usepackage[english]{babel}
\usepackage{xcolor}
\usepackage[T1]{fontenc}
\usepackage[normalem]{ulem}
\usepackage{caption}

\makeatletter
\def\paragraph{\@startsection{paragraph}{4}%
  \z@{3pt}{-\fontdimen2\font}%
  {\normalfont\bfseries}}
\makeatother

\usepackage{amssymb}
\usepackage{amsfonts}
\usepackage[tikz]{bclogo}
\usetikzlibrary{arrows, arrows.meta}
\usepackage{fullpage}
\usepackage{lmodern}
\usepackage{textcomp}
\usepackage{mathtools, bm}
\usepackage{amssymb, bm}
\usepackage{amsmath}
\usepackage[unq]{unique}
\usepackage[breaklinks]{hyperref}

\usepackage[numbers]{natbib}
\usepackage{graphicx} 

\usepackage{amsthm}
\usepackage{tikz}
\usepackage{caption}
\usepackage{subcaption}
\usepackage{cleveref}
\usepackage{thm-restate}

\newtheorem{theorem}{Theorem}[section]
\newtheorem{thm}{Theorem}
\newtheorem{lemma}[theorem]{Lemma}

\newtheorem{corollary}[theorem]{Corollary}

\newtheorem{claim}[theorem]{Claim}

\theoremstyle{definition}

\newtheorem{conjecture}[theorem]{Conjecture}
\newtheorem*{conjecture*}{Conjecture}

\newtheorem{definition}[theorem]{Definition}

\newcommand{\eps}{\varepsilon}
\newcommand{\e}{\mathrm{e}}
\newcommand{\dtv}{d_{\mathrm{TV}}}
\newcommand{\op}{\mathrm{op}}

\newcommand{\diam}{\mathrm{diam}}
\newcommand{\dist}{\mathrm{dist}}

\newcommand{\cA}{\mathcal{A}}
\newcommand{\cB}{\mathcal{B}}
\newcommand{\cC}{\mathcal{C}}
\newcommand{\cD}{\mathcal{D}}

\newcommand{\cF}{\mathcal{F}}
\newcommand{\cG}{\mathcal{G}}
\newcommand{\cH}{\mathcal{H}}

\newcommand{\cR}{\mathcal{R}}

\newcommand{\cT}{\mathcal{T}}

\newcommand{\cW}{\mathcal{W}}

\newcommand{\defined}{\mathrel{\coloneqq}}
\undef{\set}
\DeclarePairedDelimiter{\set}{\lbrace}{\rbrace}
\newcommand{\sseq}{\subseteq}
\newcommand{\inter}{\mathbin{\cap}}
\newcommand{\union}{\mathbin{\cup}}

\DeclareMathOperator{\Prob}{\mathbb{P}}

\undef{\abs}
\DeclarePairedDelimiterX{\abs}[1]
  {\lvert}{\rvert}{\ifblank{#1}{\,\cdot\,}{#1}}
\undef{\norm}
\DeclarePairedDelimiterX{\norm}[1]
  {\lVert}{\rVert}{\ifblank{#1}{\,\cdot\,}{#1}}

\title{Approximate Itai-Zehavi conjecture for random graphs}

\author[Hollom]{Lawrence Hollom}

\address{Department of Pure Mathematics and Mathematical Statistics, University of Cambridge, Cambridge, CB3 0WA, United Kingdom}
\email{lh569@cam.ac.uk}

\author[Lichev]{Lyuben Lichev}

\address{Institute of Statistics and Mathematical Methods in Economics, TU Wien, Wiedner Hauptstra\ss e 8-10, A-1040 Vienna, Austria}
\email{lyuben.lichev@tuwien.ac.at}

\author[Mond]{Adva Mond}

\address{Department of Mathematics, King’s College London, Strand, London, WC2R 2LS, United Kingdom}
\email{adva.mond@kcl.ac.uk}

\author[Portier]{Julien Portier}

\address{Ecole Polytechnique Federale de Lausanne (EPFL), CH-1015 Lausanne,
Switzerland}
\email{julien.portier@epfl.ch}

\author[Wang]{Yiting Wang}

\address{Institute of Science and Technology Austria (ISTA), 3400 Klosterneuburg, Austria}
\email{yiting.wang@ist.ac.at}

\begin{document}

\begin{abstract}
A famous conjecture by Itai and Zehavi states that, for every $d$-vertex-connected graph $G$ and every vertex $r$ in $G$, there are $d$ spanning trees of $G$ such that, for every vertex $v$ in $G\setminus \{r\}$, the paths between $r$ and $v$ in different trees are internally vertex-disjoint.
We show that with high probability the Itai-Zehavi conjecture holds asymptotically for the Erd\H{o}s-R\'enyi random graph $G(n,p)$ when $np= \omega(\log n)$ and for random regular graphs $G(n,d)$ when $d= \omega(\log n)$. 
Moreover, we essentially confirm the conjecture up to a constant factor for sparser random regular graphs.
This answers positively a question of Dragani\'{c} and Krivelevich.
Our proof makes use of recent developments on sprinkling techniques in random regular graphs.
\end{abstract}

\maketitle

\section{Introduction}

Edge-connectivity and vertex-connectivity are fundamental concepts in mathematics and theoretical computer science.
A cornerstone result in this area is Menger's theorem~\cite[Section~3.3]{Die24} which characterises $d$-(edge/vertex)-connectivity via the existence of $d$-(edge/internally vertex)-disjoint paths between every pair of vertices.
The sufficient and necessary condition given by Menger's theorem is local in the following sense: for each pair of vertices, the exhibited $d$ internally disjoint paths may have completely different shapes, and do not give much insight on the overall structure of the graph.

Beyond this local perspective, the celebrated theorem proven independently by Tutte~\cite{Tutte} and Nash-Williams~\cite{Nash-Williams} ensures that high edge-connectivity has important global consequences: every $d$-edge-connected graph contains $\lfloor d/2\rfloor$ edge-disjoint spanning trees.
In contrast, our understanding of the global structure of graphs based on their vertex-connectivity is more limited.
To address this gap, Itai and Zehavi suggested the following remarkable strengthening of Menger's theorem: 

\begin{conjecture}[Itai and Zehavi~\cite{IZ89}]\label{conj:IZ}
For every $d\ge 1$, every $d$-vertex-connected graph $G$ and every vertex $r$ in $G$, there exists $d$ spanning trees $T_1,\ldots,T_d$ of $G$ such that, for every vertex $v$ in $G\setminus \{r\}$, the paths from $v$ to $r$ in $T_1,\ldots,T_d$ are pairwise internally vertex-disjoint.
\end{conjecture}

Conjecture~\ref{conj:IZ} is immediate when $d=1$, and is known to hold for a few other small values of $d$.
Namely, the case $d=2$ was proved by Itai and Rodeh~\cite{IR88}, the case $d=3$ was confirmed independently by Cheriyan and Maheshwari~\cite{CM88}, and by Itai and Zehavi~\cite{IZ89}, and the case $d=4$ was shown to hold by Curran, Lee and Yu~\cite{CLY06}. 
A more general version of Conjecture~\ref{conj:IZ} also appears as Conjecture~6.7(b) in a paper of Frank~\cite{Fra95}. 
He also proposed an analogue for directed graphs (Conjecture~6.7(a)), which was disproved by Huck~\cite{huck1995disproof}.
\footnote{We note that, although~\cite{Fra95} and~\cite{huck1995disproof} are published within a few months from each other, at the time of Huck's work~\cite{huck1995disproof}, the conjecture of Frank had been apparently known for at least 15 years.
In fact, Huck refers a paper of Schrijver~\cite{Sch79}, who attributes this conjecture to Frank but cites a non-existent paper, obtained by a mixture of two of Frank's works. 
Since we did not find a trace of the conjecture in any of Frank's early papers, we suspect that a personal communication might have been involved.}

More generally, given a graph $G$ and a vertex $r$ in $G$, the spanning trees $T_1,\ldots,T_d$ rooted at $r$ are called \emph{independent spanning trees} (abbreviated \emph{ISTs}) if, for every vertex $v$ in $G\setminus \{r\}$, the paths from $v$ to $r$ in $T_1,\ldots,T_d$ are pairwise internally vertex-disjoint.
While $d$-regular $d$-vertex-connected graphs cannot contain $d+1$ ISTs due to simple counting restrictions (every edge may belong to at most two trees), Conjecture~\ref{conj:IZ} states that such graphs do actually contain $d$ ISTs. 
Large families of ISTs have been constructed for many particular host graphs such as hypercubes~\cite{Sha24,YTCW07}, product graphs~\cite{Che11,KWH03,OIBI96}, maximal planar graphs~\cite{NN01}, and bubble-sort networks~\cite{KPHWC19} among many others; 
for a detailed overview of the topic, see the excellent survey of Cheng, Wang and Fan~\cite{CWF23}. 

To find $\ell$ ISTs in a graph $G$, one strategy is to ensure the existence of $\ell$ disjoint connected dominating sets in $G$. 
More precisely, given disjoint connected dominating sets $D_1,\ldots,D_{\ell}$ in $G$, one can construct spanning trees $T_1,\ldots,T_{\ell}$ such that, for every $i\in \{1,\ldots,{\ell}\}$, every vertex outside $D_i$ is a leaf in $T_i$.
This way, the trees $T_1,\ldots,T_{\ell}$ rooted at an arbitrary vertex are indeed ISTs in $G$.
Currently, the best known lower bound on the number of disjoint connected dominating sets in an $n$-vertex $d$-vertex-connected graph $G$ is due to Censor-Hillel, Ghaffari, Giakkoupis, Haeupler, and Kuhn~\cite{C-HGK17} who showed that there exist at least $cd/(\log n)^2$ such sets in $G$ for some universal constant $c>0$.
Recently, Dragani\'{c} and Krivelevich~\cite{DK24} improved the latter bound to $(1+o_d(1))d/\log d$ in the case of $(n,d,\lambda)$-graphs (with $d/\lambda$ sufficiently large) and observed that this bound is optimal up to the value of the constant.
Although finding disjoint connected dominating sets allows to construct ISTs, in general, this is not the right approach towards the complete resolution of the Itai-Zehavi conjecture.
Indeed, as shown by Censor-Hillel, Ghaffari, and Kuhn \cite{censor2014new}, there exists a universal constant $C > 0$ and $d$-vertex-connected graphs on $n$ vertices which do not contain $Cd/\log n$ disjoint connected dominating sets.
Similarly, the bound proved by Dragani\'{c} and Krivelevich is tight whereas, for instance, a random $d$-regular graph is $d$-vertex-connected with high probability, see~\cite{Bol01,Wor81}.
This observation lead Dragani\'{c} and Krivelevich to pose the question of whether one could prove \Cref{conj:IZ} for random graphs or, more modestly, find $\Theta(d)$ ISTs in the random regular graph $G(n,d)$, and $\Theta(np)$ ISTs in the Erd\H{o}s-R\'enyi random graph $G(n,p)$ with $np = \omega(\log n)$.

The first main result of our work positively answers these questions.
In fact, in the case of Erd\H{o}s-R\'enyi graphs, we show that an asymptotic version of Conjecture~\ref{conj:IZ} holds. 
\begin{thm}\label{thm:main-a}
   Fix $p = p(n)\in [0,1]$ with $np = \omega(\log n)$. Then, with high probability, for every vertex $r$, the Erd\H{o}s-R\'enyi graph $G(n,p)$ contains $(1-o(1))np$ ISTs rooted at $r$.
\end{thm}

Note that \Cref{thm:main-a} is asymptotically optimal. Indeed, in the regime $np = \omega(\log n)$, routine computations show that the minimum degree (which serves as a trivial upper bound for the connectivity)
of $G(n,p)$ is equal to $(1-o(1))np$ with high probability. 
We also remark that our proof is based on the analysis of a randomised polynomial-time algorithm for finding the claimed family of ISTs.

Our second main result concerns random regular graphs.

\begin{thm}\label{thm:main-b}
    Fix $d = d(n)\in [4,n-1]$. Then, with high probability, for $(1-o(1))n$ vertices $r$, the random regular graph $G(n,d)$ contains $\lfloor d/4\rfloor$ ISTs rooted at $r$.
    Moreover, when $d = \omega(\log n)$, with high probability, for every vertex $r$, the random regular graph $G(n,d)$ contains $(1-o(1))d$ ISTs rooted~at~$r$.
\end{thm}
We make a few remarks concerning the above results. To begin with, we observe that the first part of \Cref{thm:main-b} is the first application of a recently developed sprinkling technique in random regular graphs~\cite{HLMPW25+sprinkling,IMcKSZ23} in the regime $d=\omega(1)$.

We also note that Bollob\'as~\cite{Bol01} and Wormald~\cite{Wor81} showed independently that, for every fixed $d\ge 3$, the random $d$-regular graph is $d$-vertex-connected with high probability. This result was 
extended to all $d\in [3,n]\setminus \{n-3\}$ in a sequence of works~\cite{Cooper-Frieze-Reed,Krivelevich-Sudakov-Vu-Wormald,Luc92}, with the case $d=n-3$ being exceptional due to the possibility of 4-cycles in the complement of $G(n,n-3)$.
Combining this with results from~\cite{CM88,CLY06,IZ89}, it follows that with high probability $G(n,d)$ is $d$-vertex-connected and satisfies the Itai-Zehavi conjecture for $d\in \{3,4\}$.

Finally, we observe that, with our approach, the constant $1/4$ in \Cref{thm:main-b} might be improved to $1/3$ at the price of a more technical presentation.
Since we are not able to come closer to $1$ with our approach, we stick to the current cleaner version.

\subsection{Outline of the proofs.}
The proofs for \Cref{thm:main-a} and \Cref{thm:main-b} use significantly different ideas.

\medskip

\noindent
\textbf{Erd\H{o}s--R\'{e}nyi random graphs.}
\Cref{thm:main-a} follows from the more general Theorem~\ref{thm:Gnp} which, given any $\eps > 0$, any vertex $r$ in $G\sim G(n,p)$ and $p\ge C(\log n)/n$ with suitably large constant $C = C(\eps)$, provides $(1-\eps)np$ ISTs of $G$ rooted at $r$ with probability $1-o(n^{-1})$.
For the purposes of the presentation, we assume that $p \leq \varepsilon/3$; the dense case is handled via similar but simpler methods.
The proof of Theorem~\ref{thm:Gnp} begins with a standard sprinkling argument: the graph $G(n,p)$ is decomposed into two independent random graphs $G_1\sim G(n,p_1)$ and $G_2\sim G(n,p_2)$ with $p_1$ much smaller than $p$ and $p_2\approx p$.
After fixing a root vertex $r$, we reveal its neighbours in $G = G_1 \cup G_2$ and work with $k = (1-\eps)np$ of them,
denoted by $v_1,\ldots,v_k$.
Using a breadth-first search (BFS) exploration process in $G_1$, we iteratively construct disjoint sets $C_1,\ldots,C_k$; we call them the \emph{core sets}.
Roughly speaking, in the ultimately constructed ISTs $T_1,\ldots,T_k$, the core set $C_i$ will almost coincide with the set of vertices in $T_i$ which are neither leaves nor parents of leaves.
The vertices in the core set $C_i$ will be expected to have a significant number of descendants in $T_i$ and, for this reason, will be attached as leaves in each of the trees $T_j\neq T_i$.

In more detail, we construct the core sets consecutively as follows. 
For every $i\in [k]$, upon having $C_1,\ldots,C_{i-1}$ already built, we start a BFS exploration process in $G_1\setminus (\{r,v_{i+1},\ldots,v_k\}\cup C_1\cup \dotsb\cup C_{i-1})$ from $v_i$ and run it until $\lceil\eps/3p\rceil$ vertices have been explored. 
The vertices explored during this process form our new set $C_i$.
Crucially, at the time when the BFS exploration away from $v_i$ stops, there exists a large set $B_i\subseteq C_i$ such that none of the edges between $B_i$ and $G_1\setminus (\{r,v_{i+1},\ldots,v_k\}\cup C_1\cup \dotsb\cup C_i)$ have yet been explored.

Once the core sets $C_1,\ldots,C_k$ have been constructed, the final phase of our construction consists of attaching the vertices remaining outside $(\{r\}\cup C_i)_{i=1}^k$ via paths of length at most two to the core sets.
More precisely, for every vertex $v\in V\setminus (\{r\}\cup C_i)$, if $v$ is adjacent to $B_i$ in $G$, then we incorporate $v$ in $T_i$ by connecting it to an arbitrary neighbour in $B_i$. 
If $v$ has no neighbour in $B_i$, then we use new edges exposed in $G_2$.
Namely, we look for a vertex $w_i\in V\setminus (\{r\}\cup C_i)$ of $v$ in $G_2$ which is adjacent to $B_i$. The neighbour $w_i$ then becomes the parent of $v$ in $T_i$.

Our last task is to show that, for every vertex $v\in V\setminus \{r\}$, the vertices $w_i$ required above exist and may be chosen distinct with suitably high probability.
Indeed, upon producing $k$ spanning trees, due to the disjointness of the core sets, the only possible obstruction to the independence of the said trees would be if $w_i = w_j$ for distinct indices $i,j\in [k]$.
To avoid such conflicts, we use a classic result for the existence of a perfect matching in bipartite Erd\H{o}s-R\'enyi graphs (Theorem~\ref{thm:bipartite-PM}).
Note that the analysis of this step is where the fresh randomness of $G_2$ is used; indeed, the $G_1$-neighbourhoods of many vertices in the core set $C_i$ have been fully explored but these vertices still have to be incorporated in the trees $T_j\neq T_i$.
More concretely, Theorem~\ref{thm:bipartite-PM} is applied to an auxiliary graph with parts $X\subseteq [k]$ (consisting of the indices $i$ with $N_G(v)\cap B_i=\varnothing$) and $Y=N_{G_2}(v)\setminus (\{r\}\cup C_1\cup\dotsb\cup C_k)$ where $x\in X$ is adjacent to $y\in Y$ if $y\in N_x$.

\medskip\noindent\textbf{Random regular graphs.}
\Cref{thm:main-b} is deduced from \Cref{thm:main-a} when $d = \omega(\log n)$ via a stochastic comparison result (\Cref{thm:sandwich}).
We therefore focus on the sparse regime $d \le (\log n)^2$.
Our approach here is significantly more technical and, unlike the case of Erd\H{o}s--R\'{e}nyi random graphs, this outline only aims to give a high-level idea for the proof.

The strategy from \Cref{thm:main-a} where randomness is exposed in two stages does not have a natural analogue for random regular graphs. 
To overcome this inconvenience, we first restrict our attention to even $n$ and define an alternative model obtained as a union of $d$ uniformly chosen perfect matchings on the same set of $n$ vertices conditioned on being edge-disjoint. The odd $n$ case is addressed separately; we delay the details to the end of the section.

We make use of a stochastic comparison
result (\Cref{thm:Contiguity-Disjoint-Matchings}) which says that, if a property holds with high probability in the mentioned alternative model, then it holds with high probability for the random $d$-regular graph as well.
The alternative random graph model described above provides a convenient decomposition into $d$ edge-disjoint matchings. 
We divide these matchings into $k = \lfloor d/4\rfloor$ groups of four matchings.
For every $i\in [k]$, the first three matchings in the $i$-th group are combined into a $3$-regular graph $G_i$ while the fourth matching is denoted $M_i$.
Although the 3-regular graphs $G_1,\ldots,G_k$ are neither uniform nor independent, we show an asymptotic decorrelation property which allows us to efficiently estimate the probability that $G_1,\ldots,G_k$ satisfy various useful properties, see \Cref{lem:matchings}. For example, we show that typically each graph among $G_1,\ldots,G_k$ has small diameter (\Cref{thm:BF82}).

By choosing a random vertex $r$, we construct trees $T_1,\ldots,T_k$ by a breadth-first search exploration of the graphs $G_1,\ldots,G_k$ away from $r$. These trees are typically almost the desired ISTs, meaning that there are very few vertices $v$ such that the paths from $v$ to $r$ in two distinct trees $T_i\neq T_j$ share an internal vertex.
We call such vertices \emph{bad}.
We show that, typically, all issues caused by bad vertices $v$ are resolved by exchanging the edges from $v$ and its descendants towards their parents in $T_i$ (or in $T_j$) with the edges in $M_i$ (or in $M_j$) containing the respective vertices. We call this a \emph{rerouting} operation. 

When analysing the rerouting procedure, we need to be able to compute (among other things) the probabilities of the following two types of events:
\begin{enumerate}[(i)]
    \item\label{it:1} There is a vertex $w$ appearing on the path from $v$ to $r$ in each of the trees $T_i$ and $T_j$.
    \item\label{it:2} A vertex $v$ is matched to a vertex $z$ in the matching $M_i$ (conditionally on $\bigcup_{j=1}^k G_j$ and some edges of $M_i$ being already exposed), and the paths from $z$ to $r$ in $T_i$ and $v$ to $r$ in $T_j$ share a vertex other than $r$. 
\end{enumerate}
The probabilities of these events are analysed in \Cref{cl:stable}.
To estimate the probability of these events, we introduce an idea of random overlay of the graphs $G_1,\ldots,G_k$ (also appearing in a more general form in the recent paper~\cite{HLMPW25+sprinkling} by the authors).
We call the unlabelled copy of a graph $G\subseteq K_n$ obtained by erasing the labels of $V(G)$ its \emph{skeleton}.
Given $m=O(1)$ unlabelled $3$-regular graphs $H_1,\dots, H_m$, each with $n$ vertices, we randomly overlay $H_1,\dots,H_m$ in the following way.
For each $i\in [m]$, we assign labels to the vertices of $H_i$ according to a uniformly random permutation $\sigma_i$ of the vertices of $K_n$, thus forming the labelled graph $\sigma_i(H_i)\subseteq K_n$, where the permutations $\sigma_1,\ldots,\sigma_m$ are chosen independently.
Then, for reasons related to the Poisson paradigm, $\sigma_1(H_1)\cup\dotsb\cup \sigma_m(H_m)$ does not contain any multiple edges with probability bounded away from $0$ and depending only on $m$ (\Cref{lem:disjoint}).
From here, we deduce that an event happening with high probability in this random overlay model also holds with high probability in the original model.

Then, we observe that each of the events \ref{it:1} and \ref{it:2} 
can be certified by revealing only a bounded number of graphs $G_i$ and perfect matchings $M_j$.
By combining this observation with \Cref{lem:matchings} and the overlay argument, a sequence of simple probabilistic computations and union bounds (\Cref{lem:unique,lem:low,lem:safe}) imply that, with high probability, after the rerouting procedure, neither \ref{it:1} nor \ref{it:2} occurs for any vertex of the graph.

For odd $n$ (and even $d$), we reduce the problem to the even $n-1$ case with additional constraints.
To transition between $G(n-1,d)$ and $G(n,d)$, we design an operation $\op$ which transforms a $d$-regular graph on $n-1$ vertices containing an induced matching $M$ on $d/2$ edges into a $d$-regular graph on $n$ vertices. 
More precisely, $\op$ removes the edges in $M$ from a $d$-regular graph on $n-1$ vertices $L_{n-1}$ and introduces a new vertex $v$ which connects to all vertices in $M$.
Conversely, to go from a $d$-regular graph $L_n$ on $n$ vertices to a graph on $n-1$ vertices and a matching on $d/2$ edges, one can remove a vertex from $L_n$ whose neighbourhood forms an independent set and add an arbitrary perfect matching on $N(v)$.
Equipped with this operation, we show that applying it to the random graph $G(n-1,d)$ and a random induced matching $M$ on $d/2$ edges in it produces an approximately uniformly random $d$-regular graph on $n$ vertices.
By adapting the proof of \Cref{thm:main-b} for even $n$, we construct $k$ ISTs $T_1,\ldots,T_k$ in $L_{n-1}$ which avoid the matching $M$, while also ensuring that the new vertex $v$ can be attached as a leaf to each of $T_1,\ldots,T_k$ so that the paths from $v$ to $r$ in $T_1,\ldots,T_k$ are all internally vertex-disjoint.
This finishes the proof in the case of odd $n$.

\subsection{Notation and terminology.} 
For a positive integer $n$, we denote $[n] = \{1,\ldots,n\}$.
For real numbers $a,b,c$ with $b > 0$, we write $c = a\pm b$ to say that $c\in [a-b,a+b]$. 
Floor and ceiling notation is omitted when rounding is insignificant for the argument. 
We denote the set of natural numbers by $\mathbb N = \{1,2,\dots\}$.

We work with mostly standard graph-theoretic and probabilistic notation.
For $p\in [0,1]$ and an integer $n\ge 1$, the \emph{Erd\H{o}s-R\'enyi random graph $G(n,p)$} is obtained from the complete graph on $n$ vertices by keeping edges independently with probability $p$.
The \emph{bipartite Erd\H{o}s-R\'enyi graph $G(m,n,p)$} is obtained similarly from the complete bipartite graph with two parts of sizes $m$ and $n$.
The binomial distribution with parameters $n$ and $p$ is denoted $\mathrm{Bin}(n,p)$.
We say that a sequence of events $(E_n)_{n\ge 1}$ holds \emph{with high probability} if their probability tends to 1 as $n\to \infty$.

For two probability distributions $\mu,\nu$ on (a subset of) $\mathbb R$, we say that $\mu$ \emph{stochastically dominates} $\nu$ if there is a coupling $(X,Y)$ with $X\sim \mu\text{ and }Y\sim \nu$ such that $\mathbb P(X\ge Y) = 1$.
Furthermore, for random graphs $G$ and $H$, we say that $G$ \emph{stochastically dominates} $H$ if there is a coupling of the distributions of $G$ and $H$ so that $H\sseq G$ with probability 1. For example, $G(n,p)$ stochastically dominates $G(n,p')$ whenever $p\geq p'$.

For integers $n\ge 1$ and $d\in [n-1]$, the \emph{random $d$-regular graph $G(n,d)$} is sampled uniformly from the family $\cG_d = \cG_d(n)$ of $d$-regular graphs on $n$ vertices.
For integers $k\ge 1$ and $d_1,\ldots,d_k\ge 1$ with $d_1+\dotsb+d_k\le n-1$, we denote by $G(n,d_1)\oplus \dotsb \oplus G(n,d_k)$ the union of independent copies of $G(n,d_1),\ldots,G(n,d_k)$ on the same vertex set conditionally on being edge-disjoint.

For a graph $G = (V,E)$ and a vertex set $S\subseteq V$, the \emph{neighbourhood of $S$} is the set $N(S) = \{w\in V\setminus S: w \text{ adjacent to some vertex }v\in S\}$.
We write $e(G)$ for the number of edges of $G$.

We use standard asymptotic notations.
For functions $f=f(n)$ and $g=g(n)$ (with $g$ positive), we write $f=O(g)$ to mean that there is a constant $C$ such that $|f(n)|\le Cg(n)$, and we write $f=\Omega(g)$ to mean that there is a constant $c>0$ such that $f(n)\ge cg(n)$ for sufficiently large $n$.
We write $f=\Theta(g)$ to mean that $f=O(g)$ and $f=\Omega(g)$, and we also write $f=o(g)$ to mean that $f(n)/g(n)\to 0$ as $n\to\infty$ and $f=\omega(g)$ to mean that $f(n)/g(n)\to\infty$ as $n\to\infty$.
By default, our asymptotic variable is $n$.

\subsection{Organisation of the paper.} In Section~\ref{sec:prelims}, we present some preliminary results and techniques.
In Section~\ref{sec:Gnp}, we prove \Cref{thm:main-a} along with the case of $d=\omega(\log n)$ in part~(b). Section~\ref{sec:Gnd} is then dedicated to the proof of the sparse case of \Cref{thm:main-b}.
We conclude with some open questions and possible directions for future work in \Cref{sec:conclusion}. 
Finally, \Cref{sec:appendix} contains the proof of a technical result (\Cref{thm:BF82}) used in the proof of \Cref{thm:main-b}.

\section{General preliminaries}\label{sec:prelims}

This section contains some preliminary results used in the proof of \Cref{thm:main-a} and \Cref{thm:main-b}.

\medskip
\noindent
\textbf{Counting subgraphs with fixed degrees in dense graphs.}
Our analysis of sparse random regular graphs uses the following result on subgraph enumeration in dense graphs due to McKay~\cite{McK85}.

\begin{theorem}[see Theorem~4.6 in~\cite{McK85}]\label{thm:McKay}
Fix $\varepsilon \in (0,2/3)$.
Fix a graph $X$ with vertex set $[n]$ and degree sequence $x_1,\ldots,x_n$, and let $g_1,\ldots,g_n$ be a sequence of positive numbers with an even sum. Define
\[e(G) = \frac{1}{2}\sum_{i=1}^n g_i,\quad \lambda = \frac{1}{4e(G)} \sum_{i=1}^n g_i(g_i-1)\quad \text{and}\quad \mu = \frac{1}{2e(G)} \sum_{ij\in E(X)} g_ig_j.\]
Suppose that 
\[\hat\Delta := 2+ \max_{i\in [n]} g_i\bigg(\frac{3}{2}\max_{i\in [n]} g_i + \max_{i\in [n]} x_i + 1\bigg)\le\varepsilon\sum_{i=1}^n g_i.\]
Then, the number of graphs with degree sequence $g_1,\ldots,g_n$ on the vertex set of $X$ that are edge-disjoint from $X$ is equal to
\[\frac{(2e(G))!}{e(G)! \, 2^{e(G)}\prod_{i=1}^n g_i!} \exp\left(- \lambda - \lambda^2 - \mu - O\bigg(\frac{\hat\Delta^2}{e(G)}\bigg)\right),\]
where the constant in the $O(\cdot)$ term is uniform over the choice of all parameters.
\end{theorem}

Following a notational convention from~\cite{IMcKSZ23}, for a graph $G$, we denote by $\cR_1(G)$ the family of perfect matchings in $G$.
Theorem~\ref{thm:McKay} has the following useful corollary for the count of perfect matchings in very dense regular graphs.
\begin{corollary}[Special case of Theorem~\ref{thm:McKay} for perfect matchings]\label{cor:McKay}
Fix $d = d(n)\in [3,n/5]$ and let $\bar{X}$ be any $(n-d)$-regular graph. Then,
\[|\cR_1(\bar{X})| = \frac{n!}{(n/2)! \, 2^{n/2}} \exp\left(-\frac{d-1}{2}-O\bigg(\frac{d^2}{n}\bigg)\right).\]
\end{corollary}
We note that more refined estimates for the number of perfect matchings exist (see~\cite{Cuckler-Kahn}) but~\Cref{cor:McKay} is sufficient for our purposes.

\medskip
\noindent
\textbf{Chernoff bound.} Next, we present a version of the classic Chernoff bound, see e.g.\ Theorem~2.1 in~\cite{JLR00}.

\begin{lemma}\label{lem:Chernoff}
Given a binomial random variable $X$ with mean $\mu$, for every $\delta\in (0,1]$,
\[\Prob\left(X\leq (1-\delta)\mu \right)\leq \exp\bigg(-\frac{\delta^2 \mu}{2}\bigg) \quad\text{ and }\quad \Prob\left(X\geq (1+\delta)\mu \right)\leq \exp\bigg(-\frac{\delta^2 \mu}{3}\bigg).\]
\end{lemma}

\medskip\noindent\textbf{Matchings in bipartite Erd\H{o}s-R\'enyi graphs.} The following theorem quantifies the probability that the bipartite Erd\H{o}s-R\'enyi graph does not possess a perfect matching.

\begin{theorem}[see Theorem 4.1 from \cite{JLR00}]\label{thm:bipartite-PM}
Fix $p = p(n)\in [0,1]$. Then, for every $m\le n$,
\begin{equation*}
\Prob\left(G(m,n,p)\text{ has a matching of size $m$} \right)\ge \Prob\left(G(n,n,p)\text{ has a perfect matching} \right) = 1-O(ne^{-np}).
\end{equation*}
\end{theorem}

\medskip\noindent\textbf{Comparison of random graphs.} 
The following result follows from more general theorems by Gao, Isaev and McKay~\cite[Theorem 1.1(b)(c)]{GIM-SODA}
\footnote{A more general form of \Cref{thm:sandwich} can be found in~\cite{GIM-SODA}, which is a version of the journal paper~\cite{GIM22} appearing in the Proceedings of the Fourteenth Annual ACM-SIAM Symposium on Discrete Algorithms (SODA 2020). 
However, \Cref{thm:sandwich} cannot be deduced from~\cite{GIM22} in full generality. See also the follow-up preprint by the same set of authors~\cite{GIM20+}.}
, which allows us to compare $G(n,p)$ and $G(n,d)$ for $d = \omega(\log n)$ and $p \approx d/n$:
\begin{theorem} \label{thm:sandwich}
Fix $d = d(n) = \omega(\log n)$, $\beta = \beta(n)$ tending to $0$ suitably slowly and $p = (1-\beta)d/n$. Then, there exists a coupling of $G_1\sim G(n,p)$ and $G_2\sim G(n,d)$ such that $G_1\subseteq G_2$ with high probability.
\end{theorem}

\medskip\noindent\textbf{Total variation distance.} Given two probability measures $\mu_1$ and $\mu_2$ defined on a common measurable space $(\Omega, \cF)$, the total variation distance between $\mu_1$ and $\mu_2$ is defined as  
\begin{align}\label{eq:coupl}
d_{\mathrm{TV}}(\mu_1,\mu_2)
= \sup_{A\in \mathcal F} |\mu_1(A)-\mu_2(A)|
=\inf\left\{\mathbb P(X\neq Y): (X,Y)\text{ is a coupling of }\mu_1 \text{ and }\mu_2 \right\}.
\end{align} 
Moreover, when the sample space $\Omega$ is finite, the above definition is equivalent to
\begin{align*}
d_{\mathrm{TV}}(\mu_1,\mu_2) = \frac{1}{2}\sum_{x\in \Omega} |\mu_1(x)-\mu_2(x)|.
\end{align*} 

\medskip\noindent\textbf{Contiguity.} Let $(\Omega_n, \cF_n)_{n\ge 1}$ be a sequence of measurable spaces and, for every $n\ge 1$, let $\mu_n$ and $\nu_n$ be two probability measures on the space $(\Omega_n, \cF_n)$.
We say that $(\mu_n)_{n\ge 1}$ and $(\nu_n)_{n\ge 1}$ are \emph{contiguous} if, for every sequence of events $(A_n)_{n\ge 1}$ such that $A_n\in \cF_n$ for all $n\ge 1$, we have that
\[\mu_n(A_n)\xrightarrow[n\to \infty]{} 1\iff \nu_{n}(A_{n})\xrightarrow[n \to \infty]{} 1\,.\]
For convenience, we will often abuse terminology and say that $\mu_n$ and $\nu_n$ are contiguous. 
It is easy to see that $\mu_n$ and $\nu_n$ are contiguous if $\dtv(\mu_n, \nu_n) = o(1)$.
However, the converse does not hold in general, as contiguity only concerns sequences of events which hold with high probability.
For instance, by setting $\Omega_n = \{0,1\}$, $\mu_n = \mathrm{Ber}(1/3)$ and $\nu_n = \mathrm{Ber}(2/3)$ for all $n\ge 1$, the measures $\mu_n$ and $\nu_n$ are contiguous but $\dtv(\mu_n, \nu_n) = 1/3$ for all $n\ge 1$. 

One of the key tools we will use is the following one-sided contiguity result obtained by the authors in~\cite[Theorem 1.5]{HLMPW25+sprinkling}. 
\begin{theorem}
\label{thm:Contiguity-Disjoint-Matchings}
Consider even $n$ and $d = d(n)\in [3,n^{1/10}]$.
Fix a sequence $(A_n)_{n\in 2\mathbb N}$ where $A_n$ is a set of $d$-regular graphs on $n$ vertices.
Let $\mu_d = \mu_{d,n}$ be the uniform probability distribution on $\cG_d(n)$ and let $\nu_d = \nu_{d,n}$ be the probability distribution on the same family such that, for $G\in \cG_d(n)$, $\nu_d(G)$ is proportional to the number of $1$-factorisations of $G$.
If $\lim_{n\to \infty} \nu_{d}(A_n) = 0$, then $\lim_{n\to \infty}\mu_d(A_n) = 0$. 
\end{theorem}

\medskip\noindent\textbf{Number of Triangles in $G(n,d)$.}
Finally, we make use an estimate following from the work of Gao~\cite[Theorem 9]{Gao23b}:
\begin{theorem}
\label{cla:estimate-1}
Consider integers $n$ and $d=d(n)=o(n^{2/5})$. 
Then, the number $X$ of triangles in the random $d$-regular graph $G(n,d)$ satisfies $\mathbb{E}[X]=O(d^3)$.
\end{theorem}

\medskip\noindent\textbf{Automorphism group of the random regular graphs.} Given a (labelled or unlabelled) graph $G$, an automorphism of $G$ is a permutation of vertex set $\sigma: V(G)\rightarrow V(G)$ such that $ij\in E(G)$ if and only if $\sigma(i)\sigma(j)\in E(G)$. 
The set of all automorphisms defines a group where the group operation is composition. We denote by $\mathrm{aut}(G)$ the automorphism group of $G$. 
We say the automorphism group is \emph{trivial} if the only element in $\mathrm{aut}(G)$ is the identity permutation. 
Equivalently, if $\mathrm{aut}(G)$ is trivial, then every labelling of an unlabelled copy of $G$ yields a distinct labelled graph. 
    
The last result in this section estimates the number of unlabelled $d$-regular graphs on $n$ vertices whose automorphism group is not trivial. It follows from a result due to McKay and Wormald~\cite[Corollary 3.10]{mckay1984automorphisms}.
    \begin{theorem}\label{thm:unlabelled-automorphism}
        Fix any $\varepsilon > 0$.
        For all $d=d(n)\in [3,n^{1/2-\varepsilon}]$, almost all unlabelled $d$-regular graphs on $n$ vertices have a trivial automorphism group.
    \end{theorem}

\section{\texorpdfstring{Proof of Theorem~\ref{thm:main-a}: ISTs in Erd\H{o}s-R\'enyi graphs}{Proof of Theorem~\ref{thm:main-a}}}\label{sec:Gnp}
This section is dedicated to the proof of the following statement, which we use to deduce~\Cref{thm:main-a}.

\begin{theorem}\label{thm:Gnp}
For every $\eps \in (0,1/2)$, there is a constant $C=C(\eps) > 1$ such that, for every $p = p(n) \in [C(\log n)/n,1-\varepsilon]$, the following holds with high probability: for every vertex $r$ in $G\sim G(n,p)$, $G$ contains at least $\lceil (1-\eps)np\rceil$ independent spanning trees rooted at $r$.
\end{theorem}

\begin{proof}[Proof of \Cref{thm:main-a} assuming \Cref{thm:Gnp}]
If $p = p(n)$ is bounded away from $1$, then applying \Cref{thm:Gnp} with $\eps\in (0,1-\limsup p(n))$ arbitrarily small finishes the proof. Suppose that $p = 1 - o(1)$ and fix any sequence of $(\delta_t)_{t\ge 1}$ with $\delta_t\in (0,1/2)$ for all $t\ge 1$ and $\lim_{t\to \infty} \delta_t = 0$.
Then, for every fixed $t$ and large enough $n$, one can naturally couple $G(n,1-\delta_t)$ and $G$ so that $G(n,1-\delta_t)\subseteq G$ with probability~1.
As a result, by \Cref{thm:Gnp} applied for $G(n,1-\delta_t)$ and $\eps = \delta_t$, with high probability, for every vertex $r$ in $G$, $G$ contains at least $\lceil (1-\delta_t)^2n\rceil$ ISTs rooted at $r$.
Thus, for every $\eps > 0$, by choosing $t$ suitably large, one can find $\lceil (1-\eps)np\rceil$ ISTs in $G$ as described above
with high probability, as desired. 
\end{proof}

For simplicity, we fix $G\sim G(n,p)$ with vertex set $V$ in the entire section.
We start with a standard preliminary lemma showing that $G$ has good expansion properties.

\begin{lemma}
\label{lem:expansion}
For every $c\in (0,1/2]$, there is a constant $C=C(c) > 1$ such that, for every $p = p(n)\in [C(\log n)/n,1]$, each of the following holds with probability $1-o(n^{-2})$:
\begin{enumerate}
    \item[\emph{(i)}] every set $S\sseq V(G)$ of size at most $\lceil c/p\rceil$ satisfies $|N(S)| \geq |S|(n-|S|)p/2$, and
    \item[\emph{(ii)}] every two disjoint sets of size at least $\lceil cn/7\rceil$ are connected by at least one edge.
\end{enumerate}
\end{lemma}
\begin{proof}
We start with point (i).
For every $s\in [1,\lceil c/p\rceil]$ and a vertex set $S$ of size $s$, we have 
\[
    \Prob\left(|N(S)| \leq (n-s)sp/2\right) = \Prob\left(\mathrm{Bin}(n-s, 1 - (1-p)^s)\leq (n-s)sp/2\right)
\]
Note that the inequality $(1-t)^s\le 1-ts+t^2\tbinom{s}{2}$ holding for any $t\in [0,1]$ implies that 
\[1-(1-p)^s\ge (1-(s-1)p/2)sp\ge (1-c/2)sp.\] 
Combining this with Chernoff's bound (\Cref{lem:Chernoff}) and the inequalities $s\leq n/2$ and $1-c/2\ge 3/4$ yields
\[
   \Prob\left(|N(S)| \leq (n-s)sp/2\right) \leq \Prob\left(\mathrm{Bin}(n-s, (1-c/2)sp) \leq (n-s)sp/2\right) \leq \exp\left(-\frac{snp}{72}\right).
\]

By summing over all $s\in [1,\lceil c/p\rceil]$ and using that $\tbinom{n}{s}\le (\e n/s)^s$ for all such $s$, we obtain that the said event fails with probability
\[\sum_{s=1}^{\lceil c/p\rceil} \binom{n}{s} \exp\bigg(-\frac{snp}{72}\bigg) \le \sum_{s=1}^{\lceil c/p\rceil} \exp\bigg(s\log\bigg(\frac{\e n}{s}\bigg)-\frac{snp}{72}\bigg)\le \sum_{s=1}^{\lceil c/p\rceil} \exp\bigg(-\frac{snp}{144}\bigg) = o(n^{-2}),\]
where we used that $\log(\e n/s)\le np/144$ for a suitably large $C$. This completes the proof of (i).

Moreover, a union bound shows that (ii) fails with probability at most 
\[(2^n)^2 \mathbb P(\mathrm{Bin}(c^2n^2/49, p) = 0) = 4^n (1-p)^{c^2n^2/49}\le 4^n \e^{-c^2n^2p/49} = o(n^{-2}),\]
where we used the inequality $1-p\le \e^{-p}$. This concludes the proof.
\end{proof}

We are ready to show \Cref{thm:Gnp} following the outline in the introduction.

\begin{proof}[Proof of \Cref{thm:Gnp}]
Fix $p_1 = 0.01\varepsilon^2 p/(1-p+0.01\varepsilon^2 p)$ and $p_2 = (1-0.01\varepsilon^2)p$ so that $(1-p_1)(1-p_2)=1-p$, and note that $p_1\geq 0.01\eps^2 p$.
By a standard sprinkling argument, one can decompose $G$ into independent random graphs $G_1 \sim G(n,p_1)$ and $G_2 \sim G(n,p_2)$ so that $G\sim G_1\cup G_2$.
Set $k \coloneqq \lceil (1-\eps)np\rceil$. 
As suggested in the outline, we will use $G_1$ to construct the core sets $C_1,\ldots, C_k$ and the connecting layers $N_1,\ldots, N_k$. 
Then, for every $i\in [k]$, $G_2$ will be used to construct stars centred at the vertices in the connecting layer $N_i$ and covering $V\setminus (\{r\}\cup C_i\cup N_i)$.

\medskip\noindent\textbf{Phase 1. Fixing the common root.} Fix an arbitrary vertex $r$ of $G$ before revealing any of $G_1,G_2$.
By Chernoff's bound (\Cref{lem:Chernoff}) for suitably large $C = C(\eps)$, with probability $1-o(n^{-1})$, $r$ has at least $k$ neighbours. 
We condition on this likely event, expose the neighbours of $r$ in $G$ and define $R \defined \{v_1,\dots, v_k\}$ to be a set of arbitrary $k$ neighbours of $v$.

\medskip\noindent\textbf{Phase 2. Building the core sets.} 
We construct the disjoint core sets $C_1,\dotsc,C_k$ one by one. 
Suppose that, for some $i\in [k]$, the sets $C_1,\ldots,C_{i-1}$ have been built.
Define the set
\begin{align*}
V_i\defined V\setminus (\set{r}\union R\union C_1\cup \dotsb\cup C_{i-1}).
\end{align*}
We construct the set $C_i$ via a breadth-first search (BFS) in $G_1[V_i]$, starting from $v_i$ and stopping precisely when $\lceil \eps/3p\rceil$ vertices have been added to $C_i$.
We now describe this BFS exploration in detail; note that all explored edges at the end of this step belong to $G_1[V_i]$.

\medskip\noindent\textbf{BFS algorithm.} Order the vertices in $V_i$ arbitrarily.
We maintain and gradually update three sets of vertices: the \emph{discovered set} $D_i$, the \emph{boundary set} $B_i$ and the \emph{unexplored set} $U_i$.
It is practical to think of $B_i$ as a queue processed according to the first-in-first-out order.
If $\eps\le 3p$, 
then return $D_i = \{v_i\}$, $B_i=\{v_i\}$ and $U_i = V_i\backslash \{v_i\}$. 
Otherwise, let $D_i = \{v_i\}$, $B_i = \varnothing$ and $U_i = V_i\setminus \{v_i\}$.
The BFS algorithm runs in \emph{steps} executed as long as $\abs{B_i \union D_i} <\lceil \eps/ 3p\rceil$.
At each step of the algorithm, we fix the vertex $v$ last added to $D_i$ and process it as follows.
\begin{itemize}
    \item If $v$ has unexplored edges towards $U_i$, then explore the edge $uv$ where $u\in U_i$ is the first vertex (in the order on $V_i$) for which this has not been done.
    \begin{itemize}
        \item[$\circ$] If $uv\in E(G_1)$, update \[B_i\leftarrow B_i\cup \{u\}\quad\text{and}\quad U_i\leftarrow U_i\setminus \{u\}\] 
        and proceed to the next step.
        \item[$\circ$] If $uv\not\in E(G_1)$, we simply proceed to the next step without changing any of the sets.
    \end{itemize}
    \item If all edges from $v$ to $U_i$ have been explored, look for the vertex $w$ first added to $B_i$ amongst the vertices currently present there. If $B_i=\varnothing$, we terminate the algorithm with failure. Else, we update
    \[D_i\leftarrow D_i\cup \{w\}\quad\text{and}\quad B_i\leftarrow B_i\setminus \{w\}\] 
    and proceed to the next step.
    \qed
\end{itemize}
In the sequel, we refer to the sets $B_i$, $C_i = B_i\cup D_i$ and $D_i$ produced after the termination of the BFS algorithm. 

\vspace{1em}

The next claim collects several properties of the described BFS exploration.

\begin{claim}\label{claim:Ci-property}
With probability $1-o(n^{-1})$, for every $i\in [k]$, 
\begin{itemize}
\item the produced set $C_i$ has size $\lceil \eps/3p \rceil$, and
\item $|B_i|\geq \lceil\eps/6p\rceil$.
\end{itemize}
Moreover, the following hold deterministically.
\begin{itemize}
\item $C_i\cap C_j = \varnothing$ for all $i\neq j$,
\item for all $i\in [k]$, the edges between $B_i$ and $U_i$ remain unexplored in $G_1$, and
\item for all $i\in [k]$, the BFS algorithm reveals a spanning tree $T_i'$ of $G_1[C_i]$ rooted at $v_i$ and such that all vertices in $B_i$ are leaves of $T_i'$.
\end{itemize}
\end{claim}
\begin{proof}
The last three points follow immediately from the construction.
We now prove the first two points.
If $\eps \leq 3p$, the statement holds trivially. 
We focus on the case $\eps > 3p$, implying that $p<\eps/3<1/6$.

For the first point, note that, if the BFS algorithm never returns failure, then each of the sets $C_1,\ldots,C_k$ has size exactly $\lceil \eps/3p\rceil$ by construction.
Thus, it suffices to show that the BFS algorithm does not terminate with failure with probability at least $1-o(n^{-1})$. 
For every $i\in [k]$, using that $|C_i|\leq \lceil\eps/3p\rceil$ and $p < 1/6$, we conclude that
\begin{align*}
    N \coloneqq|V_i| &\ge n-(1+k+(k-1)\lceil \eps/3p\rceil) \\
    &\ge n-(1+k+(k-1)(\eps/3p+1 )) \\
    &\ge (1- 2(1-\eps)p-(1-\eps)\eps/3)n-2 \ge  n/3.
\end{align*}
We now claim that the graph $G_1[V_i]$ is connected with probability $1-o(N^{-2}) = 1-o(n^{-2})$, which in particular shows that, with probability $1-o(n^{-1})$, the BFS algorithm does not terminate with failure.
Indeed, on the one hand, by part (i) of Lemma~\ref{lem:expansion} for $c=1/4$, any set $S$ of size $s\leq s_0 := \lceil 1/4p_1\rceil$ expands to at least $s(N-s)p_1/2 > 0$ vertices with probability $1-o(n^{-1})$, so there is no component of size at most $s_0$.
On the other hand, when $s = s_0 < \lceil N/28\rceil$, we deduce that the set $S$ expands to at least  $s_0(N-s_0)p_1/2 > \lceil N/28\rceil$ vertices, which means that with probability $1-o(n^{-1})$, connected components of size between $s_0$ and $\lceil N/28\rceil$ do not exist.
Moreover, part (ii) of Lemma~\ref{lem:expansion} implies that, with probability $1-o(n^{-1})$, there are no components of size between $\lceil N/28\rceil$ and $N-\lceil N/28\rceil$.
A union bound over the $O(n)$ choices for $i\in [k]$ finishes the proof of the first point.
        
For the second point, let $w$ be the last vertex added to $D_i$, and let $D_i' \coloneqq D_i \setminus \{w \}$.
If $|D_i'| =0$, then $|D_i| =1$, and therefore $|B_i| = \lceil \eps/3p \rceil - 1 \geq \lceil \eps/6p\rceil$, as desired.
Thus, we may assume that $|D_i'| \geq 1$.
By applying part (i) of \Cref{lem:expansion} for $G_1[V_i]$ with $c=1/2$ and choosing $C$ sufficiently large, we deduce that, with probability $1-o(n^{-2})$,
\begin{align*}
|N_{G_1[V_i]}(D_i')| \ge (N-|D'_i|)|D_i'|p_1/2 \geq |D_i| + 1.
\end{align*}
Since $N_{G_1[V_i]}(D_i') \setminus\{ w \} \sseq B_i$, we have that $\abs{B_i} \ge \abs{D_i}$. Using that $|B_i| + |D_i| = \lceil \eps/3p \rceil$, this gives $|B_i| \geq \lceil \eps/6p\rceil$, proving the second point.
\end{proof}

\medskip\noindent\textbf{Phase 3. Completing the ISTs.}
In this phase, we use as yet unexposed edges of $G$ to finish the construction of the ISTs.
In particular, we expose $G_2\sim G(n,p_2)$ with $p_2 = (1-0.01\eps^2)p$ for the first time in this process.

Recall that $k = \lceil (1-\varepsilon)np\rceil$. Then, by the Chernoff bound (\Cref{lem:Chernoff}), $r$ has at least $k$ neighbours in $G$ with probability $1-o(n^{-1})$. 
For each vertex $v\in V\setminus \{r\}$, our goal is to find distinct vertices $w_1,\ldots,w_k$ in $N_G(v)$ where, for every $i\in [k]$:
\begin{itemize}
    \item If $v\in C_i$, then $w_i$ is the parent of $v$ in $T_i'$.
    \item If $v\notin C_i$ but $N_G(v)\cap B_i\neq \varnothing$, then $w_i$ is an arbitrary neighbour of $v$ in $B_i$.
    \item If $v\notin C_i$ and $N_G(v)\cap B_i = \varnothing$, then 
    $w_i\in N_G(B_i)\cap N_{G_2}(v)\cap S$ where $S \defined V\setminus (\{r\}\cup C_1\cup \dotsb\cup C_k)$. 
\end{itemize}
Note that, by choosing $w_i$ to be the parent of $v$ in $T_i$ for every $i\in [k]$ as above, the paths from $v$ to $r$ in $T_1,\ldots,T_k$ are internally vertex-disjoint.
Indeed, for $i\neq j$, by disjointness of the core sets, the only vertices which may be shared by the paths between $v$ and $r$ in $T_i$ and in $T_j$ belong to $N_G(B_i)\cap N_G(B_j)\cap S$. In particular, this would imply that $w_i = w_j$, contradicting the fact that $w_1, \ldots, w_k$ are distinct.

Given a vertex $v \in V\setminus \{r\}$, to construct the vertices $w_1,\ldots,w_k$, we first expose the edges between $v$ and $B_1\cup \cdots\cup B_k$ in $G$. 
Denote by $I_v'\subseteq [k]$ the set of indices $i$ where $N_G(v)\cap B_i\neq \varnothing$; for every such index, $w_i$ is a neighbour of $v$ in $B_i$.
Then, define the set $I_v$ by setting $I_v=I_v'$, if $v\in S$, and otherwise set $I_v=I_v'\cup \{i\}$ where $i\in [k]$ is the unique index such that $v\in C_i$.

After exposing the edges between $v$ and $S$ in $G_2$, we consider two cases.

\medskip\noindent\textbf{Case 3.1: $\eps \ge 3p$.} We will conclude by applying \Cref{thm:bipartite-PM} to a suitable auxiliary bipartite graph.

For each $i\in [k]$, define $N_i \defined N_{G_1}(B_i)\cap S$ (see Figure~\ref{fig:proof-fig} for an illustration).

As edges between $B_i$ and $N_i$ have not yet been explored for all $i\in [k]$, and the sets $B_1,\ldots,B_k$ are pairwise disjoint,
the events $\{y\in N_i: y\in S, i\in [k]\}$ are mutually independent; indeed, they depend on disjoint sets of edges in $G_1$.
Moreover, by~\Cref{claim:Ci-property}, we have $|B_i| \ge \lceil \eps/6p \rceil$. Using this and the inequality $1-t \le e^{-t}$, we obtain that, for all $y \in S$,
\begin{align}\label{eq:ni-containment-prob}
\Prob(y\in N_i) \ge 1 - (1 - p_1)^{\abs{B_i}} \geq 1-\e^{-p_1\abs{B_i}} \geq 1-\e^{-0.01\eps^2p\cdot \lceil\eps/6p\rceil} \ge \alpha\eps^3
\end{align}
for each $i\in [k]$ and some sufficiently small absolute constant $\alpha > 0$ (e.g. $\alpha = 0.001$).
The bound \eqref{eq:ni-containment-prob} is key to our application of Theorem~\ref{thm:bipartite-PM}
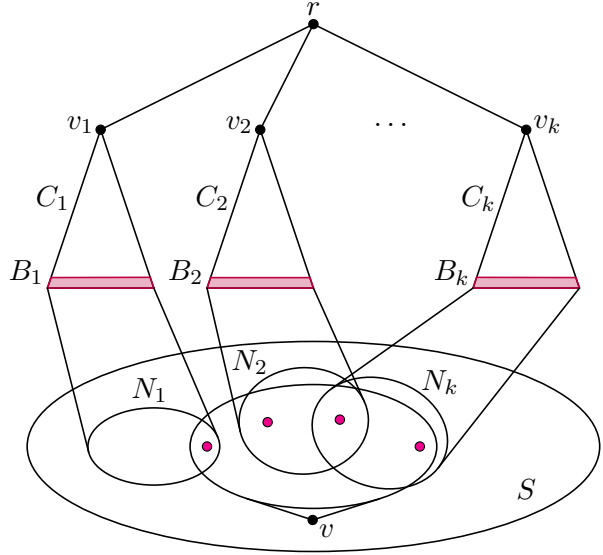
\begin{figure}[htb]
    \centering
\begin{tikzpicture}[scale=0.7,line cap=round,line join=round,x=1cm,y=1cm]
\clip(-14,-6) rectangle (11.403615967180858,4.7);
\fill[line width=0.6pt,color=purple,fill=purple,fill opacity=0.3] (-6.933403239592937,-0.8) -- (-7,-1) -- (-5,-1) -- (-5.067974144288499,-0.8) -- cycle;
\fill[line width=0.6pt,color=purple,fill=purple,fill opacity=0.3] (-3.9339292181187226,-0.8) -- (-4,-1) -- (-2,-1) -- (-2.0660922914858486,-0.8) -- cycle;
\fill[line width=0.6pt,color=purple,fill=purple,fill opacity=0.3] (1.067802634944107,-0.8) -- (1,-1) -- (3,-1) -- (2.9335327091645995,-0.8) -- cycle;
\draw [line width=0.6pt] (-6,2)-- (-7,-1);
\draw [line width=0.6pt] (-7,-1)-- (-5,-1);
\draw [line width=0.6pt] (-5,-1)-- (-6,2);
\draw [line width=0.6pt] (-3,2)-- (-4,-1);
\draw [line width=0.6pt] (-4,-1)-- (-2,-1);
\draw [line width=0.6pt] (-2,-1)-- (-3,2);
\draw [line width=0.6pt] (2,2)-- (1,-1);
\draw [line width=0.6pt] (1,-1)-- (3,-1);
\draw [line width=0.6pt] (3,-1)-- (2,2);
\draw [line width=0.6pt,color=purple] (-6.933403239592937,-0.8002097187788104)-- (-7,-1);
\draw [line width=0.6pt,color=purple] (-7,-1)-- (-5,-1);
\draw [line width=0.6pt,color=purple] (-5,-1)-- (-5.067974144288499,-0.7960775671345033);
\draw [line width=0.6pt,color=purple] (-5.067974144288499,-0.7960775671345033)-- (-6.933403239592937,-0.8002097187788104);
\draw [line width=0.6pt,color=purple] (-3.9339292181187226,-0.8017876543561675)-- (-4,-1);
\draw [line width=0.6pt,color=purple] (-4,-1)-- (-2,-1);
\draw [line width=0.6pt,color=purple] (-2,-1)-- (-2.0660922914858486,-0.8017231255424537);
\draw [line width=0.6pt,color=purple] (-2.0660922914858486,-0.8017231255424537)-- (-3.9339292181187226,-0.8017876543561675);
\draw [line width=0.6pt,color=purple] (1.067802634944107,-0.7965920951676795)-- (1,-1);
\draw [line width=0.6pt,color=purple] (1,-1)-- (3,-1);
\draw [line width=0.6pt,color=purple] (3,-1)-- (2.9335327091645995,-0.8005981274937979);
\draw [line width=0.6pt,color=purple] (2.9335327091645995,-0.8005981274937979)-- (1.067802634944107,-0.7965920951676795);
\draw [line width=0.6pt] (-6,2)-- (-2,4);
\draw [line width=0.6pt] (-2,4)-- (-3,2);
\draw [line width=0.6pt] (-2,4)-- (2,2);
\draw [rotate around={0:(-2,-4)},line width=0.6pt] (-2,-4) ellipse (5.381094528084061cm and 1.9890144092429856cm);
\draw [rotate around={0:(-2,-4)},line width=0.6pt] (-2,-4) ellipse (2.3185815104942833cm and 1.172953631140612cm);
\draw [rotate around={0:(-5,-4)},line width=0.6pt] (-5,-4) ellipse (1.2386404711415913cm and 0.7309105394984173cm);
\draw [line width=0.6pt] (-7,-1)-- (-6.225446989568821,-4.10639634893573);
\draw [line width=0.6pt] (-5,-1)-- (-3.7986840636891217,-3.8219232610833767);
\draw [rotate around={2.0213649403560248:(-2.1822444406320227,-3.5152189903082647)},line width=0.6pt] (-2.1822444406320227,-3.5152189903082647) ellipse (1.2156770702366213cm and 1.0081361059613658cm);
\draw [rotate around={-18.69667386122281:(-0.7516527345515137,-3.7456282191859724)},line width=0.6pt] (-0.7516527345515137,-3.7456282191859724) ellipse (1.2973523114632515cm and 1.0263703930968489cm);
\draw [line width=0.6pt] (-4,-1)-- (-3.376327525345874,-3.7165665136448673);
\draw [line width=0.6pt] (-2,-1)-- (-1.0413132084330374,-3.15503966857433);
\draw [line width=0.6pt] (1,-1)-- (-1.6575595084378374,-2.9037440966421078);
\draw [line width=0.6pt] (3,-1)-- (0.3453180121503498,-4.40529445889208);
\draw [line width=0.6pt] (-2.0145065770778006,-5.392285558653131)-- (-3.2594544841610857,-4.984815450600006);
\draw [line width=0.6pt] (-2.0145065770778006,-5.392285558653131)-- (-0.7611370132944938,-4.991477398132623);
\begin{scriptsize}
\draw [fill=black] (-2,4) circle (2.5pt);
\draw[color=black] (-2,4.3) node {\large{$r$}};
\draw [fill=black] (-6,2) circle (2.5pt);
\draw[color=black] (-6.3802776433301505,2.1) node {\large{$v_1$}};
\draw [fill=black] (-3,2) circle (2.5pt);
\draw[color=black] (-3.3949278405358205,2.1) node {\large{$v_2$}};
\draw[color=black] (-0.5,2.1) node {\large{$\dots$}};
\draw [fill=black] (2,2) circle (2.5pt);
\draw[color=black] (2.4,2.1) node {\large{$v_k$}};
\draw[color=black] (-6.9,0.7) node {\large{$C_1$}};
\draw[color=black] (-3.9,0.7) node {\large{$C_2$}};
\draw[color=black] (1.1,0.7) node {\large{$C_k$}};
\draw[color=black] (-7.4,-0.7) node {\large{$B_1$}};
\draw[color=black] (-4.4,-0.7) node {\large{$B_2$}};
\draw[color=black] (0.6,-0.7) node {\large{$B_k$}};
\draw [fill=magenta] (-4,-4) circle (2.5pt);
\draw [fill=magenta] (0,-4) circle (2.5pt);
\draw[color=black] (2,-4.8) node {\large{$S$}};
\draw[color=black] (-5.079910551315442,-2.9281383592724533) node {\large{$N_1$}};
\draw [fill=magenta] (-2.8611834121610187,-3.5391815422445823) circle (2.5pt);
\draw [fill=magenta] (-1.5033054691030274,-3.491256438371947) circle (2.5pt);
\draw[color=black] (-3.2,-2.4) node {\large{$N_2$}};
\draw[color=black] (0.37796822925333257,-2.7816181235524864) node {\large{$N_k$}};
\draw [fill=black] (-2.0145065770778006,-5.392285558653131) circle (2.5pt);
\draw[color=black] (-1.75,-5.6) node {\large{$v$}};
\end{scriptsize}
\end{tikzpicture}
    \caption{Illustration of the proof of Theorem~\ref{thm:main-a}. Note that $v\in S$ only for illustration convenience.}
    \label{fig:proof-fig}
\end{figure}

Next, fix a vertex $v\in V\setminus\set{r}$ and define the set $Y_v = N_{G_2}(v) \inter S$.
By \Cref{claim:Ci-property} and the assumption $3p \leq \eps$, we obtain that
\begin{equation*}
    \abs{S} =n-1-k\lceil\eps/3p\rceil \geq n-1-((1-\varepsilon)np + 1)(\eps/3p+1) 
    \geq (1-2\eps/3)n.
\end{equation*}
Thus, by choosing the constant $C$ suitably large, and by using the Chernoff bound (\Cref{lem:Chernoff}), we obtain that $|Y_v|\ge k$ with probability $1-o(n^{-2})$. We condition on this event.

Define the bipartite graph $H_v$ with parts $X = [k]\setminus I_v$ and an arbitrary subset $Y'_v \subset Y_v$ with $|Y'_v| = k$ by putting an edge between $i\in X$ and $y\in Y'_v$ if $y\in N_i$.
Then, our task is equivalent to finding a matching $M$ in $H_v$ covering $X$.
By \eqref{eq:ni-containment-prob}, $H_v$ stochastically dominates a random bipartite graph in which
every edge is present independently with probability at least $\alpha\eps^3$, where we recall that $\alpha > 0$ is an absolute constant.
In other words, $H_v$ dominates the random graph $G(m,k,\alpha \varepsilon^3)$, where $m = k-|I_v|$ and $k \ge m$.
By applying \Cref{thm:bipartite-PM}, we deduce that it contains a perfect matching with probability at least $1 - O (k\e^{-\alpha\eps^3 k})$, where the error term is of order $o(n^{-2})$ provided that $C$ is large enough in terms of $\eps$.
By a union bound over $v\in V\setminus \{r\}$, with probability $1-o(n^{-1})$, our construction finds $k$ ISTs of $G$ rooted at $r$.
Finally, a union bound over $r\in V$ finishes the proof of \Cref{thm:Gnp} in this case.

\medskip\noindent\textbf{Case 3.2: $\eps < 3p$.}
Then, $|C_1|=\cdots=|C_k|=1$, $|R| = k = \Omega(n)$ and $|S| = n-1-|R| = \Omega(n)$.
Fix any $\delta, \delta'\in (0,\eps)$ with $\delta' < \delta$.
By the Chernoff bound (\Cref{lem:Chernoff}) and a union bound, for every vertex $v\in V\setminus \{r\}$, with probability $1-o(n^{-2})$, we obtain
\begin{equation}\label{eq:ineqs} |N_{G_2}(v)\cap S| \geq (1-\delta)|S|p \qquad \text{and} \qquad |N_G(v)\cap R| = (1\pm\delta')kp. 
\end{equation}
Fix an arbitrary $v\in V\setminus \{r\}$ and recall that $I_v' \subseteq [k]$ is the set of indices $i$ where $vv_i \in E(G)$, and that $|I_v|-|I_v'|\le 1$.
We now turn to define $w_i$ for $i \in [k]\setminus I_v$.
Given the relations in~\eqref{eq:ineqs}, 
\begin{equation}\label{eq:ineq'}
\eps^3 n/6\le k - (1+\delta')kp-1\le k-|I_v| \le k - (1-\delta')kp\le (1-\delta)|S|p\le |N_{G_2}(v)\cap S|,
\end{equation}
where the first inequality used that $\delta' \le \varepsilon$ and $p \le 1-\varepsilon$, and thus $(1 - (1+\delta')p)k \ge \eps^3n/6$.
We define an auxiliary bipartite graph $H_v$ with parts $X\defined [k]\setminus I_v$ and $Y_v \subseteq N_G(v)\cap S$ of size $|Y_v| = m := k-|I_v|$ where we put an edge between $i\in I_v$ and $w\in Y_v$ whenever $wv_i\in E(G)$.
Then, $H_v \sim G(m,m,p)$.
Our goal corresponds to finding a matching in $H_v$ covering $X$. 
By~\eqref{eq:ineq'}, applying \Cref{thm:bipartite-PM}, $H_v$ contains a perfect matching with probability $1 - O(me^{-mp}) = 1-o(n^{-2})$.
A union bound over all $v\in V\setminus \{r\}$ and $r\in V$ finishes the proof for this case as well.
\end{proof}

\section{\texorpdfstring{Proof of \Cref{thm:main-b}: ISTs in random regular graphs}{Random regular graphs}}\label{sec:Gnd}
This section is dedicated to the proof of~\Cref{thm:main-b}.
We first deduce the second half where $d = \omega(\log n)$ from~\Cref{thm:main-a} and~\Cref{thm:sandwich}. 

\begin{proof}[Proof of~\Cref{thm:main-b} for $d = \omega(\log n)$]
    By~\Cref{thm:sandwich}, for some $p = (1-o(1))d/n$, there is a coupling of $G_1 \sim G(n,p)$ and $G_2 \sim G(n,d)$ such that $G_1 \subseteq G_2$ with probability $1-o(1)$.
    By~\Cref{thm:main-a}, with high probability, for every vertex $r$, there exist $(1-o(1))np = (1-o(1))d$ ISTs rooted at $r$ in $G_1$. Since $G_1 \subseteq G_2$, those trees also exist in $G_2 \sim G(n,d)$. Altogether, with probability $1-o(1)$, we have found $(1-o(1))d$ ISTs rooted at $r$, which finishes the proof.
\end{proof}

In the remainder of this section, we complete the proof of~\Cref{thm:main-b}. To this end, we show the following theorem.

\begin{theorem}\label{thm:d-reg}
Fix $d = d(n)\in [4,(\log n)^2]$ and $G = G(n,d)$. 
Then, with high probability, for all but $o(n)$ vertices $r\in V$, there are at least $\lfloor d/4\rfloor$ ISTs of $G$ rooted at $r$.
\end{theorem}

In Section~\ref{subsec:1}, we assume that $n$ is even and prove Theorem~\ref{thm:d-reg} under this additional restriction.
Then, we show the general statement in Section~\ref{subsec:2}.

\subsection{\texorpdfstring{Proving Theorem~\ref{thm:d-reg} for even $n$}{Proving Theorem~\ref{thm:d-reg} for even n}}\label{subsec:1}

By \Cref{thm:Contiguity-Disjoint-Matchings}, it is sufficient to show that, for every $d = d(n)\in [4,(\log n)^2]$, with high probability, for $n-o(n)$ vertices $r \in V(G)$, the random graph $G'(n,d)$ distributed according to $\nu_{d}$ on the space $\cG_d(n)$ contains $\lfloor d/4\rfloor$ ISTs rooted at $r$.
We first prepare the ground with several preliminary lemmas.

\subsubsection{\emph{\textbf{Some properties of random 1-factorisations}}}\label{subsubsec:usefullemmas}
For convenience of notation, we enlarge our probability space.
We define $\cG''_d(n)$ to be the family of ordered $d$-tuples of edge-disjoint perfect matchings in $K_n$.
For every $i\in [d]$ and element $H\in \cG''_d(n)$, the edges in the $i$-th matching are given \emph{colour $i$}: in particular, $H$ is identified with a properly edge-coloured $d$-regular graph.
We denote by $G''(n,d)$ a uniformly chosen edge-coloured graph in the set $\cG''_d(n)$.
Moreover, since every set of $d$ edge-disjoint matchings corresponds to $d!$ elements in $\cG''_d(n)$, one can couple the uniform probability distributions on the family of sets of $d$ edge-disjoint perfect matchings and on $\cG''_d(n)$ so that the sampled graphs coincide with probability 1. 
The next lemma is a valuable tool allowing us to compute probabilities of events concerning graphs induced by a subset of the colours.

\begin{lemma}\label{lem:matchings}
Fix $d = d(n)\in [4,n-1]$ with $d^3=o(n)$, $k\in [d]$ and colours $1\le i_1<\dotsb<i_k\le d$ inducing a $k$-regular subgraph $G^*$ from $G''(n,d)$.
Fix a sequence of events $(A_n)_{n\in 2\mathbb N}$ such that $A_n\subseteq \cG_k''(n)$ for all $n\in 2\mathbb N$.
Then, by identifying the colours $i_1,\ldots,i_k$ with $[k]$, $\mathbb P(G^*\in A_n) = (1+o(1))\mathbb P(G''(n,k)\in A_n)$.
\end{lemma}
\begin{proof}
For convenience, assume that $\{i_1,\ldots,i_k\} = [k]$. 
For every graph $H\in \cG_k''(n)$, denote by $\psi(H)$ the number of ways to extend $H$ to an element of $\cG_d''(n)$. 
For any $H\in \cG''_k(n)$, by definition, we have $\Prob(G''(n,k) = H) = 1/|\cG''_k(n)|$ and $\mathbb P(G''(n,k)\in A_n) = |A_n|/|\mathcal G''_k(n)|$. It remains to compute $\Prob(G^* \in A_n)$. Note that 
\[\mathbb P(G^* = H) = \frac{\psi(H)}{\sum_{F\in \cG_k''(n)} \psi(F)}.\]
Moreover,
by $(d-k)$ consecutive applications of \Cref{cor:McKay}, we have that, for every two graphs $H',H''\in \cG_k''(n)$, $\psi(H') = (1+O(d\cdot d^2/n))\psi(H'')$.
As a result, we deduce that 
\begin{align*}
    \mathbb P(G^*\in A_n) = \sum_{H \in A_n} \frac{\psi(H)}{\sum_{F\in \cG_k''(n)} \psi(F)}= \sum_{H \in A_n} \frac{(1+O(d^3/n))}{|\cG''(n,k)|} = (1+O(d^3/n))\mathbb P(G''(n,k)\in A_n)\,,
\end{align*}
which finishes the proof since $d^3=o(n)$.
\end{proof}

Now, we fix $k \coloneqq \lfloor d/4\rfloor$ and, for every $i\in [k]$, we denote by $G_i$ the union of the matchings in colours $4i-3,4i-2$ and $4i-1$, and denote by $M_i$ the matching in colour $4i$.
We also need the following technical lemma which computes the probability of any edge belonging to a matching in a given colour, conditionally on the edges revealed so far. 
See \cite[Theorem 1]{GO23} for a related result in the random $d$-regular graph.

\begin{lemma}
\label{lem:matching-prob}
    Fix $d = d(n)\in [4,n-1]$ with $d^3=o(n)$, $i\in [k]$ and a graph $H$ obtained by exposing $G_1,\ldots,G_k$ and $t\in [n/4]$ edges forming a matching $M'_i\subseteq M_i$.
    Then, for every $u,v\in V\setminus V(M'_i)$ with $uv\notin H$, we have 
    \[\Prob(uv\in M_i\mid H\subseteq G''(n,d)) = \frac{1 + O(d^3/n)}{n-2t}\,.\]
\end{lemma}

We remark that this result could be proved for $t\in [cn]$ for any constant $c<1/2$; we have chosen $1/4$ for simplicity.

\begin{proof}
For all $u,v\in V\setminus V(M'_i)$ with $uv\notin H$, denote by $\psi(uv)$ the number of extensions of $H\cup \{uv\}$ to an element of $\cG''_d(n)$.
Note that, by symmetry, it suffices to prove the result with $i=1$.
We apply \Cref{thm:McKay} to $H\setminus V(M_1' \union \set{uv})$ to compute the number of extensions of $H\union \set{uv}$ to $H \union M_1$, and then apply \Cref{cor:McKay} consecutively for $j\in [2,k]$ to compute the number of extensions of $H\union M_1$ to $H\union M_1 \union \dotsb \union M_k$.
Finally, a number of further applications (between 0 and 3) of \Cref{cor:McKay} allow us to compute the number of extensions of $H\union M_1 \union \dotsb \union M_k$ to an element of $\cG_d''(n)$.

By considering the error terms in the above sequence of extensions, we deduce that, for every $u,v,w\in V\setminus V(M_1')$ such that $uv, uw\notin H$, we have that $\psi(uv)/\psi(uw) = \exp(O(d^3/n)) = 1 + O(d^3/n)$.
As $u$ could match to $n-2t-O(d)$ different partners at this point, the conclusion follows.
\end{proof}
We also need a result that estimates the diameter of coloured random regular graphs of bounded degree after deletion of few random edges. This property is used only in the transition between even and odd $n$.
Its proof is inspired by the approach of Bollob\'as and Fernandez de la Vega~\cite{BF82} to analyse the diameter of the random multigraph generated from the $d$-regular configuration model, and can be found in~\Cref{sec:appendix}.

\begin{restatable}{theorem}{diameterthm}\label{thm:BF82}
Fix $\eps > 0$, $d\ge 3$ and define $s = s(n)$ to be the least integer satisfying $(d-1)^{s-3}\ge (16+\eps)dn\log n$.
    Consider a graph $H$ obtained from $G''(n,d)$ by deleting up to $(\log n)^2$ edges chosen uniformly at random.
    Then, with probability $1-o((\log n)^{-2})$, $H$ has diameter at most $s$.
\end{restatable}

\subsubsection{\emph{\textbf{The random overlay model}}}\label{subsubsec:random-overlay}
The tools developed in \Cref{subsubsec:usefullemmas}
are insufficient to estimate the probability of some important events in the model $G''(n,m)$. 
A key part of our argument consists in defining yet another random graph model, in which we can compute probabilities more efficiently:

\begin{definition}[Random overlay model]\label{def:random-overlay}
    Fix $m \in \mathbb N$ and $m$ unlabelled graphs $H_1,\dots, H_m$ on $n$ vertices.
    Consider attributing labels in $[n]$ to the vertices of each of $H_1,\ldots,H_m$ uniformly at random and independently. 
    Denote by $G$ the labelled graph with vertex set $[n]$ obtained as a union of the randomly labelled versions of $H_1,\ldots,H_m$. We denote the probability distribution of $G$ by $G^o(n, (H_1,\dots, H_m))$.
\end{definition}

The starting point of our analysis of the random overlay model is a lemma which relates it with the random $1$-factorisation model $G''(n,d)$. 
Recall that the \emph{skeleton} of a labelled graph is its unlabelled copy.

\begin{lemma}\label{lem:equiv-dist}
Fix $m\in \mathbb N$, $d_1,\ldots, d_m\in \mathbb N$ and, for every $j\in [m]$, set $D_j=\sum_{i=1}^j d_i$.
Form a random graph $G$ according to the following procedure. 
Consider a graph $G'\sim G''(n,D_m)$ and, for each $i\in [m]$, reveal the skeleton $H_i$ of the union of the perfect matchings with colours between $D_{j-1}+1$ and $D_j$. 
Conditionally on $H_1,\ldots,H_m$, define $G$ as the random graph with distribution $G^o(n, (H_1,\dots, H_m))$ conditioned on the absence of double-edges. 
Then, $G\sim G''(n,D_m)$.
\end{lemma}
\begin{proof}
It is a standard fact that, for a random variable $X$ distributed uniformly on a space $\Omega$ and a subspace $\Omega'\subseteq \Omega$, the distribution of $X$ conditionally on $X\in \Omega'$ is uniform on $\Omega'$.
Hence, conditionally on the skeletons $H_1,\ldots,H_m$, each of the distribution $G''(n,D_m)$ and the distribution $G^o(n, (H_1,\dots, H_m))$ conditioned on the absence of double-edges is uniform over the $m$-tuples of mutually edge-disjoint graphs $(G_1,\dots, G_m)$ with skeletons $H_1,\ldots,H_m$, respectively.
\end{proof}

The next lemma in this section justifies that randomly overlaying a bounded number of bounded-degree regular graphs on the same vertex set produces no multiple edges with probability bounded away from 0. 
The result was recently established in another paper of the authors~\cite[Lemma 2.10]{HLMPW25+sprinkling}.
\begin{lemma}\label{lem:disjoint}
Fix $m\ge 2$, $d_1,\ldots,d_m\ge 1$ and $m$ unlabelled regular graphs $(H_i)_{i=1}^m$ on $n$ vertices~with degrees equal to $(d_i)_{i=1}^m$, respectively.
Then, the probability that the random overlay model applied to $(H_i)_{i=1}^m$ produces $m$ pairwise edge-disjoint labelled graphs is $\exp(-D/2)+o(1)$ where $D = \sum_{1\le k < \ell \le m} d_k d_\ell$ and the $o(1)$ is uniform over all choices of $H_1,\ldots,H_m$.
\end{lemma}

Given $d_1,\dots, d_m\in \mathbb N$ such that $d_1 +\cdots + d_m = d$, a graph property $\mathcal C$ on $G''(n,d)$ is a subset of $m$-tuples of (labelled) graphs $(G_1,\ldots,G_m)$ where the graph $G_i$ consists of the matchings with indices between $d_1+\ldots+d_{i-1}+1$ and $d_1+\ldots+d_i$. 
Such $m$-tuples of graphs are called \emph{adapted} (to the degrees $d_1,\ldots,d_m$).
Note that we often abuse notation and identify the union of the graphs $G_1,\ldots,G_m$ with the $m$-tuple itself.
In this paper, we are interested in graph properties with bounded complexity.
Formally, we introduce the notion of $\ell$-certifiability.

\begin{definition}\label{def:lcertifable}
Fix positive integers $\ell\leq m$ and $d_1,\dots, d_m,d$ with $d_1 + \cdots + d_m = d$. Fix a property $\mathcal C\subseteq \cG''(n,d)$. 
We say $\mathcal C$ is an \emph{$\ell$-certifiable} property if, for every adapted $m$-tuple of graphs $G = (G_1,\dots, G_m)$ in $\cC$, there exists a set of indices $I\subseteq [m]$ of size $|I| = \ell$ with the following property: for each adapted $m$-tuple of graphs $G' = (G_1',\dots,G_m')$, if $(G_i')_{i\in I} = (G_i)_{i\in I}$, then $G'\in \mathcal C$.
\end{definition}
\noindent
Roughly speaking, a property $\cC$ is $\ell$-certifiable if the fact that $G=(G_1,\ldots,G_m)\in \cC$ is witnessed by some $\ell$ of the graphs $G_1,\ldots,G_m$.

Given a graph property, we define its unlabelled version as follows.

\begin{definition}\label{def:unlablled-dual}
    Given $m\in \mathbb N$, $d_1\dots, d_m\in \mathbb N$ with $d = d_1+\cdots +d _m$ and a graph property $\mathcal C \subseteq \cG''(n,d)$, define 
    $$\mathrm{U}(\mathcal C) = \{(H_1,\dots, H_m): (G_1,\dots, G_m)\in \mathcal C \text { and } H_i \text{ is the skeleton of } G_i\}\,.$$
\end{definition}

\subsubsection{\emph{\textbf{\texorpdfstring{Bad vertices, safe vertices and rerouting procedure}{Bad vertices, safe vertices and rerouting procedure}}}}\label{subsubsec:neven}
Recall that $k = \lfloor d/4 \rfloor$ and that, for any $i \in [k]$, we consider the graphs $G_i$ obtained as the union of the matchings in colours $4i-3, 4i-2$ and $4i-1$, and $M_i$ as the matching in colour $4i$.
To approach Theorem~\ref{thm:d-reg}, we first choose a uniformly random vertex $r$ called the \emph{root}.
Then, we construct $k$ trees $T_1,\ldots,T_k$ obtained via breadth-first search (BFS) exploration of $G_1,\ldots,G_k$, respectively, implemented according to an arbitrary ordering of the $n$ vertices in each graph. Note that the latter orderings are often left implicit.
While these trees are not quite ISTs in $G''(n,d)$, we will later show that they typically are such after a few modifications, using edges from $M_1, \ldots, M_k$. 
Note that given $G_1,\dots, G_k$ with corresponding vertex orderings, the trees $T_1,\dots, T_k$ are uniquely defined.
We start with a definition which plays a central role in our proof.
\begin{definition}
    For the spanning trees $T_1, \dots, T_k$ on $V$ rooted at $r\in V$ and distinct integers $i,j\in [k]$, a vertex $v\in V\setminus \{r\}$ is called \emph{$(i,j)$-bad} if there is a vertex $u\in V\setminus \{v,r\}$ such that $u$ belongs to the path from $v$ to $r$ in each of $T_i$ and $T_j$.
    In this case, the vertex $u$ will be called an \emph{anchor} for $v$.
    A vertex is called \emph{bad} if it is $(i,j)$-bad for some distinct $i,j\in [k]$, and \emph{good} otherwise.
\end{definition}
In particular, the fact that $T_1,\ldots,T_k$ are ISTs in $G''(n,d)$ is equivalent to the absence of bad vertices.
For a vertex $v$ and a tree $T_i$, we denote by $T_i(v)$ the subtree of $T_i$ formed by $v$ and the descendants of $v$, and denote by $P_i(v)$ the path from $v$ to $r$ in $T_i$.

Next, we describe a \emph{rerouting procedure} using edges from $M_1, \ldots, M_k$ to modify the trees $T_1, \ldots, T_k$ and turn them into ISTs.
These modifications are done locally to vertices which are either bad or descendants of bad vertices (as the latter may become bad after the rerouting).
For a vertex $v\in V$ and a tree $T_i$, if the path from $v$ to $r$ in $T_i$ is modified throughout the process, we say that the vertex $v$ is \emph{rerouted} in $T_i$.
For a graph $G$ and vertices $u,v$ therein, we denote by $\dist(u,v;G)$ the graph distance between $u$ and $v$ in $G$.

\medskip\noindent\textbf{Rerouting procedure.} 
Expose the graphs $G_1,\ldots,G_k$ and define $\sigma$ to be an arbitrary ordering of the bad vertices. Moreover, to every vertex $u\in V$, we assign an initially empty index set $I(u)$, the set of indices of trees in which $u$ will be rerouted. 
We process the bad vertices one by one according to $\sigma$.
Upon reaching a bad vertex $v$, for every pair of indices $i,j\in [k]$ such that $v$ is $(i,j)$-bad and for every anchor $u$ for $v$ in $T_i,T_j$, compare $\dist(u,r;T_i)$ and $\dist(u,r;T_j)$. Denote by $\ell(i,j,u)$ the index for which the distance is larger; in case of equality, fix $\ell = \min\{i,j\}$.
Then, for every vertex $w$ in $T_\ell(v)$, add index $\ell$ to $I(w)$.
Upon processing all bad vertices, for every vertex $w$ and every index $\ell\in I(w)$, replace the edge from $w$ to its parent in $T_{\ell}$ by the edge from the matching $M_\ell$ containing $w$.

\vspace{0.5em}
Note that if a vertex $u$ is such that the paths $(P_i(u))_{i=1}^k$ remain unchanged, then $u$ was good before the rerouting procedure and remains good thereafter.
However, a vertex $u$ can be good and still some of the paths $(P_i(u))_{i=1}^{k}$ may be changed during the procedure.
Indeed, this can happen if the good vertex $u$ is a descendant of an $(i,j)$-bad vertex $w$ in $T_i$ or in $T_j$.
We next define \emph{safe} vertices to be those which are never affected by reroutings as above.

\begin{definition}
We call a vertex $u$ \emph{safe} if $I(v)$ is empty, and \emph{unsafe} otherwise.
\end{definition}

The importance of safe vertices, as already mentioned, is that their paths to $r$ in $T_1,\ldots,T_k$ do not change during the rerouting procedure.
Note that the set of unsafe vertices contains the set of bad vertices and descendants of bad vertices in all trees.
In the sequel, we show that typically almost all vertices are safe.
Moreover, we show that, by rerouting all bad vertices and descendants of bad vertices in the appropriate trees, all bad vertices become good and no new bad vertices are created, thereby yielding a family of $k$ many ISTs.

\medskip\noindent\textbf{Analysing the rerouting procedure.}
Next, we prove several auxiliary lemmas concerning the rerouting procedure. We will make frequent use of~\Cref{lem:equiv-dist} and compute probabilities of various events in the random overlay model without explicit mention of \Cref{lem:equiv-dist}.

Define 
\begin{equation}\label{eq:def_psi_beta}
\psi = \psi(n) := 20\log_2\log n\qquad \text{and}\qquad \beta = \beta(n) := \log_2 n + \psi.
\end{equation}
We note that the use of base-2 logarithms in the definitions of $\psi$ and $\beta$ is due to the fact that, locally, $G(n,3)$ resembles a binary tree.

The first lemma shows that, for the trees $T_1,\ldots,T_k$ constructed via a BFS in $G_1,\ldots,G_k$, respectively, every bad vertex is $(i,j)$-bad with anchor $u$ for a unique choice of indices $i,j$ and a vertex $u$.
\begin{lemma}\label{lem:unique}
With high probability, for every bad vertex $v$, there is a unique choice of a vertex $u\in V\setminus \{v,r\}$ and a pair of distinct indices $i,j\in [k]$ such that $v$ is $(i,j)$-bad with anchor $u$.
\end{lemma}
\begin{proof}
First of all, by combining \Cref{lem:matchings} and \Cref{thm:BF82}, a union bound over the graphs $G_1,\ldots,G_k$ where $k = O((\log n)^2)$ shows that, with high probability, each of them has diameter at most $\beta$.
We condition on this event in the remaining considerations.

First, we compute the probability that there is a bad vertex $v$ with anchor $u\in V\setminus \{v,r\}$ which precedes $v$ in three of the graphs $G_1,\ldots,G_k$: a property which is 3-certifiable.
Fix distinct indices $i,j,\ell\in [k]$.
Note that equipping the skeletons of $G_i,G_j,G_\ell$ with arbitrary vertex orderings and revealing the position of the root vertex $r$ in each of them determines the skeletons of the rooted trees $T_i,T_j,T_\ell$.
Fix a vertex $v\in V\setminus \{r\}$ and expose its position in the three graphs.
Then, one can identify the paths $P_i(v),P_j(v),P_\ell(v)$ without exposing the internal vertex labels.
Recalling that each of $T_i,T_j,T_\ell$ has height at most $\beta$, and so each of the paths $P_i(v),P_j(v),P_\ell(v)$ has length at most $\beta$. 
Hence, the probability that $u$ belongs to all three of them is at most $(\beta/(n-2))^3 = n^{-3+o(1)}$. 
A union bound over $O(k^3)$ choices of indices $i,j,\ell$, $O(n)$ choices for $u$ and $O(n)$ choices for $v$ proves that this case happens with probability at most $o(1)$. 
In particular, with high probability, for any choice of index sets $\{i,j\}\neq \{s,t\}$, no vertex is $(i,j)$-bad and $(s,t)$-bad with the same anchor. 

Next, we compute the probability that there is a bad vertex $v$ with distinct anchors $u_1,u_2\in V\setminus \{v,r\}$, each preceding $v$ in two of the trees $T_1,\ldots,T_k$, say $T_i,T_j$ for $u_1$ and $T_s,T_t$ for $u_2$ with $i\neq j$ and $s\neq t$. 
Note that this is a 4-certifiable property.
Similarly to the previous case, irrespectively of possible coincidences of $i,j$ and $s,t$ (as $u_1$ and $u_2$ are distinct), the probability of the above event is $n^{-4+o(1)}$.
A union bound over $O(k^4)$ choices of indices $i,j,s,t$, $O(n^2)$ choices for $u_1,u_2$ and $O(n)$ choices for $v$ proves that this case happens with probability at most $o(1)$. 
Thus, with high probability, no vertex is bad with two distinct anchors. Together with the previous case this finishes the proof.
\end{proof}

Before we turn to the proof of \Cref{thm:d-reg}, we show two simple but useful lemmas.
Recall $\beta$ defined in~\eqref{eq:def_psi_beta}.

\begin{lemma}\label{lem:low}
With high probability, for every vertex $u\in V\setminus \{r\}$, there is at most one $i\in [k]$ such that $\dist(u,r;T_i)\le \beta/3$.
\end{lemma}
\begin{proof}
Similarly to the proof of~\Cref{lem:unique}, note that, by \Cref{lem:matchings}, \Cref{thm:BF82} and by taking a union bound over all $k = O((\log n)^2)$ graphs, with high probability, each of $G_1,\ldots,G_k$ has diameter at most $\beta$.
In particular, with high probability, each of $T_1, \ldots, T_k$ is of height at most $\beta$.
Since these trees are uniquely determined by the underlying graphs $G_1,\ldots,G_k$ equipped with arbitrary vertex orderings, respectively, 
\[\{u\in V:\quad \exists i,j\in [k]\quad\text{such that}\quad i\neq j,\quad \dist(u,r;T_i)\le \beta/3\quad\text{and}\quad\dist(u,r;T_j)\le \beta/3\}\] 
is a 2-certifiable property.
Moreover, for distinct $i,j\in [k]$ and $u\in V\setminus \{r\}$, the probability that $\max\{\dist(u,r;T_i), \dist(u,r;T_j)\}\le \beta/3$ 
is at most $(3\cdot 2^{\beta/3-1})^2/n^2 = o(1/(k^2n))$.
The conclusion follows by a union bound over the $O(k^2)$ choices for $i,j$ and $O(n)$ choices for $u$.
\end{proof}

The next lemma concerns the number and the positions of the bad vertices in the trees $T_1,\ldots,T_k$. Recall $\psi$ defined in~\eqref{eq:def_psi_beta}. 

\begin{lemma}\label{lem:safe}
With high probability, there are at most $\beta^{20}$ unsafe vertices. Moreover, with high probability, each of them is at distance at least $\beta-2\psi$ from the root $r$ in each of the trees $T_1,\ldots,T_k$.
\end{lemma}
\begin{proof}
Again, similarly to the proof of~\Cref{lem:unique}, note that, by \Cref{lem:matchings}, \Cref{thm:BF82} and by taking a union bound over all $k = O((\log n)^2)$ graphs, with high probability, each of $G_1,\ldots,G_k$ has diameter at most $\beta$ and, therefore, each of $T_1, \ldots, T_k$ has height at most $\beta$. We assume this event in the sequel.
Next, note that the property that a vertex $u$ is unsafe is 2-certifiable: indeed, it is sufficient to find indices $i,j \in [k]$ and a vertex $v$ such that $v$ is $(i,j)$-bad with $\ell(i,j,v)=i$ and $v \in P_i(u)$. 
For a fixed pair of vertices $u,v$ and indices $i,j$ as above, expose the skeletons of $G_i,G_j$ and the position of the root vertex $r$ in them, thus determining $T_i,T_j$.
Then, consecutively expose the positions of $u$ in $G_i$ and of $v$ in $G_i,G_j$.
As each of $P_i(u)$, $P_i(v)$ and $P_j(v)$ has length at most $\beta$, it follows that the probability of the event $\{v\in P_i(u)\}\cap \{v\text{ is }(i,j)\text{-bad}\}$ is bounded from above by $\beta^3/n^2$. 
A union bound over the $O(k^2)$ choices of indices and $O(n)$ choices of a vertex $v$ shows that the vertex $u$ is unsafe with probability $O(k^2\beta^3/n)$.
As a result, Markov's inequality and the fact that $k^2\beta^3 = o(\beta^{20})$ imply the first statement.

For the second statement, note that the property that there exists an unsafe vertex which is at distance at most $\beta-2\psi$ from the root $r$ in some of the trees is 3-certifiable: 
indeed, on top of the described certificate $i,j\in [k]$ ensuring that $u$ is an unsafe vertex, one needs to verify that $u$ is at distance less than $\beta-2\psi$ in some of the trees $T_1,\ldots,T_k$. 
A similar computation shows that the expected number of unsafe vertices at distance at most $\beta-2\psi$ from $r$ in some of the trees is at most $(2^{\beta-2\psi}/n)\cdot k^3\beta^3\le 2^{-\psi} \beta^{14} = o(1)$. An application of Markov's inequality finishes the proof. 
\end{proof}

\subsubsection{\emph{\textbf{Proof of \Cref{thm:d-reg} for even $n$}}}
We are now ready to prove~\Cref{thm:d-reg} for the case where $n$ is even.
\begin{proof}[Proof of \Cref{thm:d-reg} for even $n$]
We show that \Cref{thm:d-reg} holds with high probability for a uniformly randomly chosen vertex $r$; it then follows immediately that the statement holds for almost all vertices.
Assume the properties from Lemmas~\ref{lem:unique},~\ref{lem:low} and \ref{lem:safe}.
We start with the following claim regarding rerouted vertices.

\begin{claim}\label{cl:stable}
With high probability, for every vertex $w$ and every $i\in I(w)$, the vertex $z$ with $wz\in M_i$ is safe. In particular, $i\notin I(z)$.
\end{claim}
\begin{proof}
Expose the graphs $G_1,\ldots,G_k$.
On the one hand, by \Cref{lem:safe}, the number of unsafe vertices is at most $\beta^{20}$ (which, in particular, also dominates the total number of vertices rerouted in some of the trees).
On the other hand, by \Cref{lem:matching-prob}, for all $i\in [k]$, $w$ with $i\in I(w)$ and a unsafe vertex $z$, the probability that $wz\in M_i$ is $O(1/n)$. 
A union bound over $k$ choices for the index $i$ and $O(\beta^{20})$ 
choices for each of the vertices $w,z$ with $i\in I(w)$ implies the claim. 
\end{proof}
\Cref{cl:stable} shows that, with high probability, for every vertex $w$ rerouted in a tree $T_i$ to a vertex $z$ (where $wz\in M_i$), none of the paths $(P_j(z))_{j=1}^k$ is modified by the rerouting procedure. We assume this property in the sequel.

Suppose that the vertex $w$ becomes (or remains) bad after the rerouting procedure. 
In particular, this means that one of the following happens:
\begin{enumerate}[\upshape{\textbf{A\arabic*}}]
    \item\label{item:A1} there are distinct indices $\ell_1,\ell_2\in I(w)$ and vertices $z_1,z_2\in V$ such that $wz_1\in M_{\ell_1}$, $wz_2\in M_{\ell_2}$ and the paths $P_{\ell_1}(z_1),P_{\ell_2}(z_2)$ share a vertex different from $r$, or
    \item\label{item:A2} there are indices $\ell_1\in I(w)$ and $\ell_2\in [k]$
    and a vertex $z\in V$ such that $wz\in M_{\ell_1}$ and the paths $P_{\ell_1}(z),P_{\ell_2}(w)$ share a vertex different from $r$.
\end{enumerate}

Observe property \ref{item:A1} is 6-certifiable: indeed, for indices $i_1,i_2,j_1,j_2$ with $i_1\neq j_1$ and $i_2\neq j_2$, $\ell_1 \in \{i_1,j_1\}$ and $\ell_2 \in \{i_2,j_2\}$, exposing the graphs $G_{i_1},G_{j_1},G_{i_2},G_{j_2}$ and the matchings $M_{\ell_1},M_{\ell_2}$ allows us to find a vertex $v_1$ which is $(i_1,j_1)$-bad and $w\in V(T_{\ell_1}(v_1))$, a vertex $v_2$ which is $(i_2,j_2)$-bad and $w\in V(T_{\ell_2}(v_2))$, and vertices $z_1,z_2$ such that $wz_1\in M_{\ell_1}$, $wz_2\in M_{\ell_2}$ and $V(P_{\ell_1}(z_1))\cap V(P_{\ell_2}(z_2))\neq \varnothing$. 
Property \ref{item:A2} is 4-certifiable for similar reasons.

First, we show that, with high probability, \ref{item:A1} does not hold for any vertex of the graph.
Fix indices $i_1,j_1,i_2,j_2$, vertex $w\in V$, vertices $v_1$ and $v_2$ (which we will soon assume to be $(i_1,j_1)$-bad and $(i_2,j_2)$-bad, respectively), and indices $\ell_1 = \ell(i_1,j_1,v_1)$ and $\ell_2 = \ell(i_2,j_2,v_2)$. 
We consider three cases. 

\vspace{0.5em}
\noindent
\textbf{Case 1.} $\{i_1,j_1\}\cap\{i_2,j_2\} = \varnothing$. Then, the probability that $v_1$ is $(i_1,j_1)$-bad and $v_2$ is $(i_2,j_2)$-bad is $O((\beta^2/n)^2) = n^{-2 + o(1)}$. 

First, suppose that $w\notin \{v_1,v_2\}$. Then, by \Cref{lem:safe}, each of the trees $T_{\ell_1}(v_1)$ and $T_{\ell_2}(v_2)$ contains vertex $w$ independently with probability $O(2^{2\psi}/n) = n^{-1+o(1)}$.
In total, the probability that $v_1$ is $(i_1,j_1)$-bad, $v_2$ is $(i_2,j_2)$-bad and $w$ is a descendant of $v_1$ in $T_{\ell_1}$ and of $v_2$ in $T_{\ell_2}$ is at most $n^{-4+o(1)}$, and a union bound over $O(k^4)$ choices of indices and $O(n^3)$ choices for the vertices $v_1,v_2,w$ shows that with high probability, vertices with these properties do not exist.

Now, suppose that $w\in \{v_1,v_2\}$, say $w=v_1$ without loss of generality. Then, the probability that $v_1$ is a descendant of $v_2$ in $T_{\ell_2}$ is $O(2^{2\psi}/n)$ and, therefore, the probability that $v_1$ is $(i_1,j_1)$-bad, $v_2$ is $(i_2,j_2)$-bad and $v_1$ is a descendant of $v_2$ in $T_{\ell_2}$ is at most $n^{-3+o(1)}$.
A union bound over $O(k^4)$ choices of indices and $O(n^2)$ choices for the vertices $v_1,v_2$ shows that with high probability, vertices with these properties do not exist.

\vspace{0.5em}
\noindent
\textbf{Case 2.} $i_1 = i_2 = i$ and $j_1 \neq j_2$. If $w\notin \{ v_1,v_2 \}$, then, by exposing the vertices preceding $v_1$ and $v_2$ in $T_i$ first, one can see that the probability that $v_1$ is $(i_1,j_1)$-bad and $v_2$ is $(i_2,j_2)$-bad is still $O((\beta^2/n)^2)$. 
Since $\ell_1\neq \ell_2$, once again, each of the trees $T_{\ell_1}(v_1)$ and $T_{\ell_2}(v_2)$ contains vertex $w$ independently with probability $O(2^{2\psi}/n)$.
Hence, we conclude by the same union bound as in Case 1 when $w\notin \{ v_1,v_2 \}$.

Now, suppose that $w\in \{v_1,v_2\}$, say $w=v_1$ without loss of generality. Then, by consecutively checking whether $v_2$ belongs to some of the paths $P_i(v_1)$ and $P_{j_1}(v_1)$, that $P_{j_1}(v_1)$ shares a vertex with $P_i(v_1)$ and that $P_{j_2}(v_2)$ shares a vertex with $P_i(v_2)$, the total error probability is again $n^{-3+o(1)}$. Then, the conclusion follows by a union bound over $O(k^3)$ choices for the indices and $O(n^2)$ choices for the vertices $v_1,v_2$.

\vspace{0.5em}
\noindent
\textbf{Case 3.} $i_1 = i_2 = i$ and $j_1 = j_2 = j$. Note that, in this case, vertices $v_1$ and $v_2$ must be $(i,j)$-bad for distinct anchors $u_1$ and $u_2$: indeed, if the two anchors coincide, by definition of the rerouting procedure, $\ell_1$ and $\ell_2$ would also coincide. 
Due to this observation, a similar computation using \Cref{lem:safe} shows that vertices with these properties do not exist, and ensures that~\ref{item:A1} holds with high probability.

\vspace{0.5em}

Next, we show that, with high probability, \ref{item:A2} does not hold for any of the vertices of the graph.
Fix distinct indices $i,j$, a vertex $w\in V$, a vertex $v$ (which we will shortly assume to be $(i,j)$-bad), and indices $\ell_1 = \ell(i,j,v)$ and $\ell_2$.
We consider two cases.

\vspace{0.5em}
\noindent
\textbf{Case 1.} $\ell_2\notin \{i,j\}$. Then, the probability that $v$ is $(i,j)$-bad is $O(\beta^2/n)$.

Suppose that $w\neq v$. Then, by \Cref{lem:safe}, the probability that $w$ is a descendant of $v$ in $T_{\ell_1}$ is $O(2^{2\psi}/n)$. Moreover, expose the edge $wz\in M_{\ell_1}$ and the vertex labels on the path $P_{\ell_1}(z)$.
Then, the probability that $P_{\ell_2}(w)$ contains at least one of these labels is $O(\beta^2/n)$.
The total error probability being $n^{-3+o(1)}$, the conclusion follows by a union bound over $O(k^3)$ choices for the indices and $O(n^2)$ choices for the vertices $w,v$. 

Suppose that $w = v$. Again, expose $z$ and the vertex labels on the path $P_{\ell_1}(z)$.
Then, the probability that $P_{\ell_2}(v)$ contains at least one of these labels is $O(\beta^2/n)$. 
This time, the total error probability being $n^{-2+o(1)}$, the conclusion follows by a union bound over $O(k^3)$ choices for the indices and $O(n)$ choices for the vertex $v$.

\vspace{0.5em}
\noindent
\textbf{Case 2.} $\ell_2\in \{i,j\}$. Again, the probability that $v$ is $(i,j)$-bad is $O(\beta^2/n)$.
Upon this event, denote by $u$ the anchor of $v$; recall that this anchor is unique by \Cref{lem:unique}. 

Suppose that $w\neq v$. Then, by \Cref{lem:safe}, the probability that $w$ is a descendant of $v$ in $T_{\ell_1}$ is $O(2^{2\psi}/n)$.
Expose the vertex labels on the path $P_{\ell_2}(w)$.
As $\dist(u,r,T_{\ell_1}) \ge \dist(u,r,T_{\ell_2})$ by \Cref{lem:low}, we have $\dist(u,r,T_{\ell_1})\ge \beta/3$ and, therefore, the probability that $z$ is a descendant of $u$ in $T_{\ell_1}$ is $O(2^{2\beta/3}/n) = n^{-1/3+o(1)}$.
At the same time, by our assumption of \Cref{lem:unique}, $P_{\ell_1}(u)$ and $P_{\ell_2}(w)$ do not share internal vertices.
Thus, exposing vertex $z$ in $T_{\ell_1}$ and assuming that $z\notin T_{\ell_1}(u)$, the probability that $P_{\ell_1}(z)$ has a common vertex with $P_{\ell_2}(w)$ is $O(\beta^2/n)$.
As a result, the total error probability being $n^{-7/3+o(1)}$, the conclusion follows by a union bound over $O(k^3)$ choices for the indices and $O(n^2)$ choices for the vertices $w,v$. 

Finally, suppose that $w = v$ and expose vertex $z$ in $T_{\ell_1}$.
The probability that $z$ is a descendant of $u$ in $T_{\ell_1}$ is $O(2^{2\beta/3}/n) = n^{-1/3+o(1)}$; assume that this is not the case.
Then, the probability that $P_{\ell_1}(z)$ has a common vertex with $P_{\ell_2}(v)$ remains $O(\beta^2/n)$.
As a result, the total error probability is $n^{-4/3+o(1)}$, and the conclusion follows by a union bound over $O(k^3)$ choices for the indices and $O(n)$ choices for the vertices $v$, which completes the proof.
\end{proof}

\subsection{\texorpdfstring{Proving Theorem~\ref{thm:d-reg} for odd $n$}{Proving Theorem~\ref{thm:d-reg} for odd n}}\label{subsec:2}
The goal of this section is to prove \Cref{thm:d-reg} for odd $n$ and even $d\le (\log n)^2$.
We will reduce the problem to the case of even $n$ with additional constraints.

\medskip\noindent\textbf{Step 1: Finding strong ISTs in $G(n-1,d)$.}
In~\Cref{subsec:1}, we have proven that with high probability, for almost all vertices $r$ in $G(n-1,d)$, there are $k = \lfloor d/4\rfloor$ ISTs rooted at $r$.
Our first step in the analysis for odd $n$ is to strengthen the definition of ISTs for later usage.
Similarly to~\Cref{subsec:1}, by~\Cref{thm:Contiguity-Disjoint-Matchings}, we will work with $G''(n-1,d)$ instead.

For a matching $M\subseteq K_n$ and a graph $G$, we say that $M$ is an \emph{induced $d$-matching} in $G$ if it is an induced matching on $d$ vertices in $G$.  
We now introduce the concept of \emph{strong ISTs}. 
Recall that, for a rooted spanning tree $T_i$ and a vertex $v$, the path from $v$ to the root $r$ is denoted by $P_i(v)$.
\begin{definition}
Consider odd $n\to\infty$ and even $d = d(n)\in [3,(\log n)^2]$.
Let $L_{n-1}\in \cG_d''(n-1)$, $r\in V(G)$ and $M$ be an induced $d$-matching in $L_{n-1}$.
We say that $L_{n-1}$ admits $k$ \emph{strong ISTs} rooted at $r$ with respect to $M$ if there are ISTs $T_1,\ldots,T_k \in L_{n-1}$ rooted at $r$ that additionally satisfy the following two constraints:
\begin{enumerate}[\upshape{\textbf{P\arabic*}}]
    \item\label{item:P1} None of the edges of the trees $T_1,\ldots,T_k$ is in $M$.
    \item\label{item:P2} For every distinct $i,j \in [k]$ and (not necessarily distinct) $v,w \in V(M)$, the intersection of $V(P_i(v))$ and $V(P_j(w))$ is equal to $\{v,r\}$ if $w = v$, and to $\{r\}$ if $w\neq v$. 
\end{enumerate}
\end{definition}
We need the following simple fact:
\begin{lemma}\label{prop:induced-d-matching}
Consider odd $n\to\infty$ and even $d = d(n)\in [3,(\log n)^2]$.
Fix $L_{n-1}\in \cG_d''(n-1)$ and sample $d/2$ edges uniformly at random with repetition and independently from $E(L_{n-1})$. 
Then the probability that these edges do not form an induced $d$-matching in $L_{n-1}$ is at most $d^4/n$.
\end{lemma}
\begin{proof}
Denote by $S$ a uniformly random multiset of $d/2$ edges in $L_{n-1}$.
Then, it suffices to bound from above the probability that there exists a vertex such that the graph induced by the vertices at distance at most two from $v$ contains at least two edges in $S$ (counted with multiplicity).
As every such neighbourhood contains at most $d(d-1)$ edges, by a union bound, this probability is at most $n\binom{d/2}{2}(\frac{d(d-1)}{nd/2})^2 \le d^4/n$, which finishes the proof.
\end{proof}

For $\eps\in (0,1)$, we define the property
\begin{equation}\label{eqn:strong-ISTs}
\begin{split}
\cC_\eps^s = \{&L_{n-1}\in \cG_d''(n-1): \text{there are at least $(1-\eps)n$ vertices $r\in V(L_{n-1})$ and at least}\\
&\text{$(1-\eps) \tfrac{(nd/2)^{d/2}}{(d/2)!}$ induced matchings $M$ such that $L_{n-1}$ contains $k$ strong ISTs with respect to $M$} \\
&\text{and rooted at $r$ .}\}
\end{split}
\end{equation}
Next, we show that, for every $\eps > 0$, $G''(n-1,d)$ satisfies $\cC_\eps^s$ with high probability.

\begin{lemma}\label{lem:odd-n-exists-strong-ISTs}
Fix any $\eps\in (0,1)$ and $d=d(n) \in [4,(\log n)^2]$.
Then, $\Prob(G''(n-1,d) \in \cC_\eps^s) = 1-o(1)$.
\end{lemma}
\begin{proof}
Fix $L_{n-1}\sim G''_d(n-1,d)$. Sample a root vertex $r$ uniformly at random, and sample a multiset $S$ of $d/2$ edges uniformly at random with repetition from $E(L_{n-1})$. Note that, by \Cref{prop:induced-d-matching}, with high probability $S$ is an induced $d$-matching.
We are going to show that, with high probability, $G''(n-1,d)$ contains $k$ strong ISTs rooted at $r$ with respect to $S$, which is sufficient to conclude.
Recall the coloured $3$-regular graphs $G_1,\ldots,G_k$ and matchings $M_1,\ldots,M_k$ from \Cref{subsubsec:neven}: here, we define these graphs with respect to $L_{n-1}$.
We begin with the following claim.
\begin{claim}\label{cla:random-edges}
With high probability (with respect to the randomness of $S$), $S$ forms an induced $d$-matching and, for every $i\in [k]$, $|G_i\cap S| \leq 20\log \log n$ and $|M_i\cap S| \leq 10\log\log n$.
\end{claim}
\begin{proof}
By \Cref{prop:induced-d-matching}, for any outcome of $L_{n-1}$, $S$ forms an induced $d$-matching with probability at least $1-d^4/n = 1-o(1)$. 

Recall that the graphs $G_1,\dots, G_k,M_1,\dots, M_k$ are edge-disjoint.
Since every edge in $S$ is chosen uniformly at randomly and independently, for every $i\in [k]$, (the distribution of) $|G_i\cap S|$ is stochastically dominated by $\mathrm{Bin}(d/2, 3/d)$ and (the distribution of) $|M_i\cap S|$ is stochastically dominated by $\mathrm{Bin}(d/2, 1/d)$. 
By Chernoff's bound (see \Cref{lem:Chernoff}) and union bound, it follows that $|G_i\cap S| \leq 20\log \log n$ and $|M_i\cap S| \leq 10\log\log n$ for all $i\in [k]$ with probability at least $1 - (\log n)^{-5}$ for all $i\in [k]$, which finishes the proof.
\end{proof}
From now on, we fix an arbitrary outcome of $S$ satisfying \Cref{cla:random-edges} and work only with the randomness of $L_{n-1}$.
Recall from~\eqref{eq:def_psi_beta} that $\psi = 20\log_2\log n$ and $\beta = \log_2 n + \psi$.
\begin{claim}\label{cla:diameter}
With high probability (with respect to the randomness of $L_{n-1}$), for all $i\in [k]$, the diameter of $G_i\setminus S$ is at most~$\beta$.
\end{claim}
\begin{proof}
Assume that the properties from \Cref{cla:random-edges} hold. 
Condition on $m_i = |G_i\cap S|\in [0, \psi]$ for every $i\in [k]$. 
Then, $G_i\cap S$ is a uniformly chosen subset of $E(G_i)$ of size $m_i$. 
Using \Cref{lem:matchings}, it suffices to show the statement for $(G_i''\backslash S)_{i=1}^k$ where, for every $i\in [k]$, $G_i''$ is uniformly distributed over $\cG_3''(n-1)$.
Therefore, by \Cref{thm:BF82}, each $G_i''\backslash S$ has diameter at most $\beta$ with probability $1 - o(1/(\log n)^2)$. The conclusion follows by a union bound over all $i\in [k]$. 
\end{proof}

From now on, we assume that the properties from \Cref{cla:random-edges} and \Cref{cla:diameter} hold. 
We follow the proof in \Cref{subsubsec:neven} for $G_1\backslash S,\dots, G_k\backslash S$. 
The proofs of \Cref{lem:unique}, \Cref{lem:low} and \Cref{lem:safe} carry through to the odd $n$ setting without modification: indeed, they only use that the relevant graphs have maximum degree $3$ and diameter at most $\beta$. We assume those results in the sequel. 
It remains to take into account the following two additional restrictions:
\begin{enumerate}
\item The edges in $S$ cannot be used in the rerouting procedure, as per \ref{item:P1}.
\item We need to ensure that \ref{item:P2} is satisfied.
\end{enumerate}

For the first point, we know from \Cref{lem:safe} that, for each $i\in [k]$, at most $\beta^{20}$ vertices in $G_i$ need to be rerouted using $M_i$. 
Thus, at most $\beta^{20}$ edges in $M_i$ can be potentially affected. 
However, by \Cref{cla:random-edges}, for every $i\in [k]$, $|M_i\cap S|\le 10\log\log n$.
Thus, a union bound over $i\in [k]$ and at most $\beta^{20}$ endangered edges in each matching show that the first point fails with probability at most
\[O\bigg(k\cdot \beta^{20}\cdot \frac{\log\log n}{n}\bigg) = o(1).\]

For the second point, note that, for any two (possibly coinciding) vertices $u,v\in V(S)$, the property that $P_i(u)$ and $P_j(v)$ share a vertex different from $r$ (if $u\neq v$) or outside $\{v,r\}$ (if $u=v$) for some $i\neq j$ is $2$-certifiable.
By reusing the random overlay idea from \Cref{subsubsec:neven} and a union bound, the probability that there exist two vertices in $V(S)$ violating this property is at most $O(d^2\cdot k^2\cdot \beta^2/n )= o(1)$, where $d^2$ stands for the number of choices of vertices $u,v$ in $V(S)$, $k^2$ dominates the number of choices of distinct indices $i,j\in [k]$ and $O(\beta^2/n)$ dominates the probability some vertex appears on each of $P_i(u) \cup P_i(y)$ and $P_j(v) \cup P_j(z)$, where $uy \in M_i$ and $vz \in M_j$ (thus taking the possibility of rerouting into account).   

Putting everything together, for a uniformly chosen rooted vertex $r$, with high probability, $G''(n-1,d)$ admits $k$ strong ISTs with respect to a uniformly chosen multiset of $d/2$ edges, as desired.
\end{proof}

\medskip\noindent\textbf{Step 2: Transition to an unlabelled model.}
We introduce the following unlabelled model:
\begin{definition}
Fix $d\geq 3$ and denote by $\cG^u_d(n)$ the set of unlabelled $d$-regular graphs on $n$ vertices.
We define $G^u(n,d)$ to be the uniform distribution over $\cG^u_d(n)$.
\end{definition}
We say a property $\cC$ is \emph{label-independent} if the skeleton of each graph $G$ determines if $G\in \cC$.
Observe that, for every $\eps \in (0,1)$, the property $\cC_\eps^s$ is label-independent.
Recall the definition of the unlabelled version $\mathrm{U}(\cC)$ of a given property $\cC$ (\Cref{def:unlablled-dual}).

The goal of this subsection is to show the following lemma which translates any with high probability result for label-independent properties from labelled to unlabelled graphs.

\begin{lemma}\label{lem:odd-n-equiv-models}
Fix $d=d(n)\in [3,(\log n)^2]$ and a label-independent property $\mathcal C\subseteq \cG_d''(n-1)$.
We have~that $\Prob(G''(n-1,d)\in \mathcal C) = 1-o(1)$ if and only if $\Prob(G^u(n-1,d)\in \mathrm{U}(\cC)) = 1-o(1)$. 
Moreover, the same result holds for $\cG_d''(n)$, $G''(n,d)$, and $G^u(n,d)$.
\end{lemma}
\begin{proof}
By~\Cref{thm:unlabelled-automorphism}, there are $o(|\cG^u_d(n-1)|)$ unlabelled graphs with a non-trivial automorphism group. Hence, for every label-independent property $\cC$,
\[(|\mathrm{U}(\cC)|-o(|\cG^u_d(n-1)|))n!\le |\cC|\le |\mathrm{U}(\cC)|n!.\]
By using that $|\cG_d''(n-1)| = (1-o(1))|\cG^u_d(n-1)|\cdot n!$ (again implied by~\Cref{thm:unlabelled-automorphism}) and dividing the latter chain of inequalities by $|\cG_d''(n-1)|$, we obtain that
\[(1-o(1))\Prob(G^u(n-1,d)\in \mathrm{U}(\cC))\le \Prob(G''(n-1,d)\in \mathcal C)\le (1+o(1))\Prob(G^u(n-1,d)\in \mathrm{U}(\cC)),\]
yielding the first statement. The second statement follows analogously.
\end{proof}

\medskip\noindent\textbf{Step 3: The even-odd transition.}
Now, we introduce an operation $\op$ which allows to transition between $\cG^u_d(n-1)$ and $\cG^u_d(n)$.

\begin{definition}
Let $L_\varnothing$ be the unlabelled graph on $n$ vertices without edges.
Fix $L_{n-1}\in \cG^u_{d}(n-1)$ and a multiset $S$ of $d/2$ edges in $L_{n-1}$ (in particular, some edges may appear more than once in $S$).
The operation $\op$ takes $(L_{n-1},S)$ as an input and outputs a graph $L_n\in \cG^u_d(n)\cup \{L_\varnothing\}$ constructed as follows: 
\begin{itemize}
\item If $S$ is an induced $d$-matching in $L_{n-1}$, then $L_n$ is obtained by deleting $S$ from $L_{n-1}$, adding a new vertex $v$ and connecting it by edges to each of the vertices in $V(S)$. In particular, $L_n\in \cG^u_d(n)$.
\item Otherwise, $L_n = L_\varnothing$.
\end{itemize}
\end{definition}

Recall the definition of $\cC_\eps^s$ from~\eqref{eqn:strong-ISTs} and, for every $\eps\in (0,1)$, define the property $\cD_\eps$ to be the set of graphs $L_n\in \cG_d''(n)$ such that, for at least $(1-\eps)n$ vertices $r$ in $L_n$, $L_n$ contains $k$ ISTs rooted at $r$.
The goal of this part is to show the following lemma.

\begin{lemma}\label{lem:odd-even-transition}
For every $\eps > 0$, $\Prob(G^u(n,d)\in \mathrm{U}(\mathcal D_\eps))\ge 1-\eps$.
\end{lemma}
We need a version Strassen's theorem, which allows us to establish a coupling between two measures conveniently. See \cite[Corollary 2.2]{IMcKSZ23} or \cite[Proposition 6]{Kop24}.
\begin{theorem}\label{thm:strassen}
Let $\delta,\varepsilon\in [0,1]$ and $G$ be a bipartite graph with parts $S,T$.
Suppose that $S$ contains at least $(1-\delta)|S|$ vertices of degree at least $(1-\eps)e(G)/|S|$, and similarly $T$ contains at least $(1-\delta)|T|$ vertices of degree at least $(1-\eps)e(G)/|T|$.
Then, there is a coupling $(X,Z)$ with $X$ uniformly distributed on $S$ and $Z$ uniformly distributed on $T$ such that $\mathbb P(XZ\notin E(G))\le 2\delta+\eps/(1-\eps)$.
\end{theorem}
We are now ready to prove~\Cref{lem:odd-even-transition}:
\begin{proof}[Proof of~\Cref{lem:odd-even-transition}]
Define $\cH_d = \cH_d(n)$ to be the family of pairs $(L_{n-1}, S)$ with $L_{n-1}\in \cG^u_d(n-1)$ and $S$ as a multiset of $d/2$ edges in $L_{n-1}$. 
Let $\lambda$ be a measure supported on $\cH_d$ obtained by first sampling a graph $L_{n-1}\sim G^u_d(n-1)$ and then sampling a multiset $S$ consisting of $d/2$ edges in $E(L_{n-1})$ uniformly at random (with repetition).
Define $\eta_d$ to be the measure defined as $\eta_d(A) = \lambda(\op^{-1}(A))$ for all $A\subseteq \cG^u_d(n)\cup \{L_\varnothing\}$.
Let $\mu_d$ be the uniform measure on $\cG^u_d(n)$.
We extend the domain of $\mu_d$ by setting $\mu_d(L_\varnothing)=0$.
We now relate $\eta_d$ with $\mu_d$. 

\begin{claim}\label{cla:construct}
Consider $n\to\infty$ odd and even $d = d(n)\in [3,(\log n)^2]$. Then, $\dtv(\mu_d,\eta_d) = o(1)$.
\end{claim}
\begin{proof}
First of all, denote by $\mu_d'$ the uniform measure on the set $\cG^u_d(n)\cup \{L_\varnothing\}$.
It is evident that $\dtv(\mu_d,\mu_d')=o(1)$.
Moreover, since every graph $L_{n-1}\in \cG^u_d(n-1)$ participates in the same number of pairs in $\cH_d$, sampling a random $d$-regular graph on $n-1$ vertices first and a multiset of $d/2$ edges after that is equivalent to sampling a uniformly chosen element of $\cH_d$.

We will show that $d_{\mathrm{TV}}(\mu_d', \eta_d) = o(1)$. This is clearly enough to establish the desired result, as a simple application of the triangle inequality would then yield $\dtv(\mu_d,\eta_d) \leq \dtv(\mu_d,\mu_d') + \dtv(\mu_d',\eta_d) =o(1)$.
To this end, recall from~\eqref{eq:coupl} that it suffices to find a coupling of the two random variables drawn from the respective measures so that the probability they differ is~$o(1)$.

Consider the bipartite graph $\cW$ with parts $\cA=\cG^u_d(n)\cup\{L_\varnothing\}$ and $\cB=\op(\cH_d)$ with multiplicity (that is, a graph $H\in \op(\cH_d)$ has $|\op^{-1}(H)|$ copies in $\cB$). 
We add an edge between $G\in \cA$ and $H\in \cB$ if they correspond to the same graph.
    
By construction, every vertex $H\in \cB$ has degree $1$ in $\mathcal W$.
Moreover, the degree of a vertex $G\in \cA$ in $\cW$ is equal to $|\op^{-1}(G)|$. 
Conveniently, for every graph $L_n\in \cG^u_d(n)$, this is equal to the number of vertices in $L_n$ whose neighbourhood forms an independent set multiplied by the number of perfect matchings on $d$ vertices, the latter of which is $m \coloneqq d!/((d/2)!2^{d/2})$.    
By \Cref{cla:estimate-1}, we know that the expected number of triangles in a random $d$-regular graph is at most $Cd^3$ for some absolute constant $C > 0$. 
Thus, by Markov's inequality, there are at least $(1-1/(2\log n))(|\cA|-1) \geq (1-1/\log n)|\mathcal A|$ graphs in $\cA$ containing at most $2Cd^3\log n$ triangles. 
Note that, in a graph with at most $\ell$ triangles, there are at least $n - 3\ell$ vertices with neighbourhood forming an independent set. 
In particular, there are at least $(1-1/\log n)|\cA|$ graphs $G\in \cA$ with at least $(1-6Cd^3\log n/n)n\cdot m$ neighbours in $\cB$.

Finally, note that every vertex $G\in \cA$ has degree at most $mn$ in $\cW$ and, in particular, $e(\cW) / \abs{\cA} \leq mn$.
We may thus deduce that there are at least $(1-1/\log n)|\cA|$ graphs $G\in \cG^u_d(n)$ with at least $(1-6Cd^3\log n/n)e(\cW)/|\cA|$ neighbours in $\cB$.
Applying \Cref{thm:strassen} with $\eps = 6Cd^3\log n/n$ and $\delta = 1/\log n$ shows that $d_\mathrm{TV}(\mu'_d, \eta_d) = o(1)$, as desired.
\end{proof}

Next, we make an easy observation which relates strong ISTs with ISTs.
\begin{claim}\label{claim:reduce}
Consider $n\to \infty$ odd, $d = d(n)\in [3,(\log n)^2]$ even and $(L_{n-1},M)\in \op^{-1}(\cG^u_d(n))$. 
For every $t\geq 1$, if $L_{n-1}$ admits $t$ strong ISTs with respect to $M$ rooted at a vertex $r$, then $\op(L_{n-1}, M)$ also admits $t$ ISTs rooted at $r$.
\end{claim}
\begin{proof}
Define $H = \op(L_{n-1}, M)$ and let $\cT$ be a set of strong ISTs in $L_{n-1}$ with respect to $M$.
Define a set of spanning trees $\cT'$ in $H$ as follows: add a new vertex $v$ to $L_{n-1}$ and, for each tree $T\in \cT$, select a distinct matching edge $uw\in M$ and connect $v$ to either $u$ or $w$ arbitrarily. 

As the trees in $\cT'$ extend the trees in $\cT$, it is immediate that, for any $u\neq v$, the paths from $u$ to $r$ in the $k$ distinct trees in $\cT'$ are internally disjoint.
Moreover, the paths from $v$ to $r$ in the trees in $\cT'$ are internally disjoint by \ref{item:P2}, which finishes the proof. 
\end{proof}

We are ready to finish the proof of \Cref{lem:odd-even-transition}. Denote by $\mathcal H'\subseteq \mathcal H_d$ the set of pairs $(L_{n-1},S)$ such that, for at least $(1-\eps)n$ vertices $r$ in $L_{n-1}$, there are $k$ strong ISTs in $L_{n-1}$ with respect to $S$ and rooted at $r$.
By definition of $\cC_\eps^s$,
every $L_{n-1}\in \mathrm{U}(\cC_\eps^s)\subseteq \cG^u_d(n-1)$ participates in at least $(1-\eps) (nd/2)^{d/2}/(d/2)!$ pairs in $\cH'$.
By combining this with \Cref{lem:odd-n-exists-strong-ISTs} and \Cref{lem:odd-n-equiv-models} (yielding $\mathbb P(G^u(n-1,d) \in \mathrm{U}(\cC_\eps^s))=1-o(1)$), we obtain that $\lambda(\mathcal H') \ge 1-\eps-o(1)\ge 1-2\eps$.
Observe also that, for each $(L_{n-1}, S)\in \mathcal H'$, by~\Cref{claim:reduce}, there are $k$ ISTs in $\op(L_{n-1}, S)$.
As a result, we obtain that 
\[
    \mathcal H' \subseteq \{\op^{-1}(H): H\in \mathrm{U}(\cD_\eps)\}\,,
\]
implying that $\lambda(\{\op^{-1}(H): H\in \mathrm{U}(\cD_\eps)\}) \ge 1-2\eps$ and thus $\eta_d(\mathrm{U}(\cD_\eps)) \ge 1-2\eps$.
By~\Cref{cla:construct}, we deduce that $\mu_d(\mathrm{U}(\cD_\eps)) \ge 1-3\eps$. Since the latter inequality holds for every $\eps > 0$, the lemma follows.
\end{proof}

\medskip\noindent\textbf{Step 4: Finishing the proof.} 
\Cref{thm:d-reg} for odd $n$ follows by combining \Cref{lem:odd-n-equiv-models} and \Cref{lem:odd-even-transition}.

\section{Concluding remarks}
\label{sec:conclusion}
We have established the Itai-Zehavi conjecture holds asymptotically for Erd\H{o}s-R\'enyi graphs random graphs and dense random regular graphs. 
For sparse random regular graphs, we prove the conjecture approximately. 
Each of the results leaves some natural open questions.

For Erd\H{o}s-R\'enyi graphs, our result (\Cref{thm:main-a}) applies in the regime of $p = \omega(\log n/n)$. 
One may wonder if the Itai-Zehavi conjecture still holds for sparser graphs. 
Moreover, our result is only asymptotically tight: the number of ISTs we found is not equal to the minimum degree. In fact, our proof uses the fact that there is a gap between these two quantities in a significant way. 

Given the above considerations, we propose the following conjecture: 
\begin{conjecture}\label{conj:1}
Fix any $p \geq \log n/n$ and let $G\sim G(n,p)$ have minimum degree $\delta(G)$. Then, with high probability, for any vertex $r\in V(G)$, the graph $G$ contains $\delta(G)$ ISTs rooted at $r$.
\end{conjecture} 

A classic result by Bollob{\'a}s and Thomason~\cite{Bollobas-Thomason} states that the vertex connectivity of $G(n,p)$ is the same as the minimum degree for all $p\in (0,1)$. 
Therefore, the Itai-Zehavi conjecture would imply the existence of such a family of ISTs rooted at any given vertex in $G$. 
If the answer to \Cref{conj:1} is positive, techniques of Krivelevich and Samotij~\cite{krivelevich2012optimal} might be useful for its resolution.

For sparse random $d$-regular graphs, the best one can get using our method is $\lfloor d/3\rfloor$; indeed, recall that the graphs $G_1,\ldots,G_k$ explored in the proof are 3-regular. Moreover, our proof would not work with 2-regular graphs as these are typically not even connected. 
At the same time, by a simple double-counting argument, proving the existence of more than $\lfloor d/2\rfloor$ ISTs in $G(n,d)$ would require using some edges $uv$ in two trees $T_i,T_j$ in the family where $u$ is a parent of $v$ in $T_i$ and $v$ is a parent of $u$ in $T_j$. Due to this dependence, $\lfloor d/2\rfloor$ seems to be a natural barrier and proving the existence of more than $\lfloor d/2\rfloor$ ISTs in $G(n,d)$ is of interest. 

More ambitiously, a natural strengthening of this question is to find $d$ or, at least, $d-o(d)$ ISTs in a random $d$-regular graphs.
\begin{conjecture}\label{conj:2}
Fix $d = d(n)\in [3,n-1]$ and $G\sim G(n,d)$. Then, with high probability, $G$ contains $d$ ISTs if $d\neq n-3$ and $d-1$ ISTs otherwise.
\end{conjecture}

The case $d=n-3$ in \Cref{conj:2} is exceptional because, contrary to every other $d\in [3,n-1]$, the random graph $G(n,n-3)$ is not $(n-3)$-connected with probability bounded away from 0; see the discussion in the paragraph after \Cref{thm:main-b}. 
In particular, \Cref{conj:2} is a weaker version of the Itai-Zehavi conjecture for random regular graphs.

\paragraph{Acknowledgements.} 
Hollom was supported by the Internal Graduate Studentship of Trinity College, Cambridge. Lichev was supported by the Austrian Science Fund (FWF) grant No.~10.55776/ESP624. Mond was supported by UK Research and Innovation grant MR/W007320/2.  Wang was supported by the ERC Starting Grant ``RANDSTRUCT'' No.~101076777. For open access purposes, the authors have applied a CC BY public copyright license to any author accepted manuscript version arising from this submission.
Part of this research was done during a visit of the fourth author to IST Austria. We thank IST Austria for its hospitality.

\bibliographystyle{amsplain_initials_nobysame_nomr}
\providecommand{\bysame}{\leavevmode\hbox to3em{\hrulefill}\thinspace}
\providecommand{\MR}{\relax\ifhmode\unskip\space\fi MR }
\providecommand{\MRhref}[2]{%
  \href{http://www.ams.org/mathscinet-getitem?mr=#1}{#2}
}
\providecommand{\href}[2]{#2}

\appendix

\section{\texorpdfstring{Proof of \Cref{thm:BF82}}{Proof of Theorem \ref{thm:BF82}}}\label{sec:appendix}
Recall that $G''(n,d)$ denotes the random graph on $n$ vertices obtained as a union of an $d$-tuple of perfect matchings sampled uniformly at random conditionally on being edge-disjoint. 
In particular, $G''(n,d)$ is an $d$-regular graph.
For the convenience of the reader, we restate the theorem.

\diameterthm*

\begin{proof}
We prove the statement for a graph $H$ sampled according to a configuration model where every vertex is adjacent to a single half-edge in each of $d$ colours, and half-edges are matched uniformly at random conditionally on respecting the colours. 
While multiple edges in different colours are possible,
by~\Cref{lem:disjoint}, $d$ matchings chosen independently and uniformly at random are edge-disjoint with probability bounded away from 0. 
Hence, any statement holding with probability $1-o((\log n)^{-2})$ in this configuration model also holds with probability $1-o((\log n)^{-2})$ in the original model. 

First, we show that the graph sampled according to this configuration model has good expansion properties.
Then, we use this expansion to prove that the diameter of the graph does not decrease too much after deleting up to $(\log n)^2$ of its edges uniformly at random.

The first part of the proof follows the lines of the proof of~\cite[Theorem 1]{BF82}, with certain modifications.
We generate the random graph $G''(n,d)$ on the vertex set $V = [n]$ according to the following process.
Starting with a vertex $v\in V$, we expose the edges of the graph in a BFS manner.
At each step, we reveal an edge in a predetermined colour which is not yet incident to the currently processed vertex. The BFS processes each vertex exactly once, thus gradually revealing balls with growing radii around the vertex $v$ in the coloured configuration model.

Our process runs in stages.
We start by fixing the sets
\[S_0 = B_0 = \{v\},\; E_0 = \varnothing,\; R_0 = V\setminus\{v\},\]
and, for every $\ell \in [d]$, we further fix $V_0(\ell) = \varnothing$.
Fix $i\ge 0$ and assume that stage $i$ was completed.
Stage $i+1$ is described as follows.
We will consider the following sets defined with respect to the underlying graph exposed at the end of stage $i+1$, ignoring the colours:
\begin{align*}
    S_{i+1} &= \text{ the set of vertices at distance at most $i+1$ from $v$,}\\ 
    B_{i+1} &= \text{ the set of vertices at distance precisely $i+1$ from $v$,}\\
    E_{i+1} &= \text{ the set of exposed edges,}\\
    R_{i+1} &= \text{ the set of vertices not adjacent to any edge in $E_{i+1}$.}
\end{align*}
Moreover, we define the graph $H_{i+1} \defined (V,E_{i+1})$.
Lastly, consider the set of colours $[d]$. For each $\ell \in [d]$ and stage $i$, let 
\[
    V_{i+1}(\ell) = \text{ the set of vertices covered by the edges of colour $\ell$ in $H_{i+1}$.}
\]

Given $i \ge 0$ and assuming we have defined the sets $S_i, B_i, E_i, R_i$ and $V_i(\ell)$ for all $\ell \in [d]$, we describe stage $i+1$. It consists of several steps.
For every $X\in \{S,E,R\}$, set $X_{i+1}^{0} = X_i$ and further set $B_{i+1}^0 = \varnothing$.
Moreover, denote $V_{i+1}^{0}(\ell) = V_i(\ell)$ for all $\ell \in [d]$.
Assume that, for some $j \ge 0$, we have defined $X_{i+1}^{j}$ for every $X\in \{S,B,E,R\}$ and $V_{i+1}^{j}(\ell)$ for all $\ell \in [d]$.
If every $u \in B_i$ is adjacent to $d$ edges in pairwise different colours in $E_{i+1}^{j}$, then, for every $X\in \{S,B,E,R\}$, set $X_{i+1} = X_{i+1}^{j}$, $V_{i+1}(\ell) = V_{i+1}^j(\ell)$ for every $\ell \in [d]$ and terminate stage $i+1$.
Otherwise, let $u\in B_i$ be the vertex of degree smaller than $d$ and minimal label in the so-far exposed graph $H_{i+1}^{j} \defined (V, E_{i+1}^{j})$.
For a vertex $x \in V$ and $i,j\ge 0$, define $N_{i}^{j}(x)$ to be the (open) neighbourhood of $x$ in $H_{i}^{j}$, $\deg_{i}^{j}(x) \defined |N_{i}^{j}(x)|$ to be its degree, and $C_i^{j}(x)$ to be the set of colours $\ell \in [d]$ of edges incident to $x$ in $H_i^{j}$, that is, such that $x \in V_i^{j}(\ell)$.
Fix the minimal $r\in [d]\setminus C_{i+1}^{j}(u)$ and reveal the neighbour $w$ of $u$ in the matching of colour $r$: this vertex is sampled uniformly at random from the set $V \setminus (V_{i+1}^{j}(r) \cup \{u\})$. 
After adding the $r$-coloured edge $uw$ to the graph $H_{i+1}^{j}$, we update
\begin{align*}
S_{i+1}^{j+1} = S_{i+1}^{j} \cup \{w\},\quad E_{i+1}^{j+1} = E_{i+1}^{j} \cup \{uw\},\quad R_{i+1}^{j+1} = R_{i+1}^{j}\setminus \{w\},\quad V_{i+1}^{j+1}(r) = V_{i+1}^{j}(r) \cup \{u,w\}.
\end{align*}
Moreover, for all $\ell \in [d]\setminus \{r\}$, set $V_{i+1}^{j+1}(\ell) = V_{i+1}^{j}(\ell)$.
If $w \notin R_{i+1}^{j}$, set $B_{i+1}^{j+1} = B_{i+1}^{j}$ and call the edge $uw$ \emph{dispensable}.
Otherwise, if $w \in R_{i+1}^{j}$, set $B_{i+1}^{j+1} = B_{i+1}^{j} \cup \{w\}$ and call the edge $uw$ \emph{indispensable}.
We show that, with high probability, this exploration process ``expands'' rather quickly.

Recall that, by definition, the $(i+1)$-st stage of this process lasts as long as there are vertices $x \in B_i$ adjacent to less than $d$ edges in $H_{i+1}^{j}$, that is, with $\deg_{i+1}^{j}(x) < d$.
In particular, if the edge $uw$ chosen at step $j$ of stage $i+1$ is dispensable, then we have $w \in B_i\cup B_{i+1}^{j}$.

Consider the $j$-th step in the $(i+1)$-st stage of the above process, and assume that so far the process executed $k$ steps in total.
Let $u \in B_i$, $r \notin C_{i+1}^{j}(u)$ and $w \in V\setminus (V_{i+1}^{j}(r) 
\cup \{u\})$
be the chosen vertices and colour for this step.
Since there have been $k$ steps in the process so far, we have that $|B_i\cup B_{i+1}^{j}| \le k$.
Moreover, after $k$ steps, we have $|V_{i+1}^{j}(r)| \le 2k$ and thus $|V\setminus (V_{i+1}^{j}(r) \cup \{u\})| \ge n-2k-1$.
Hence, the probability that the edge $uw$ is dispensable is at most 
\begin{equation}
\label{eq:appendix:Pdis}
    \frac{|B_i\cup B_{i+1}^{j}|}{|V\setminus (V_{i+1}^{j}(r)\cup \{u\} ) |} \le \frac{2k}{n-2k-1}.
\end{equation}
Therefore, the probability that more than one dispensable edge is revealed during the first $k = O(1)$ steps is at most
\begin{equation}
\label{eq:appendix:1dis}
    \binom{k}{2} \left(\frac{2k}{n-2k-1} \right)^2 = O\left(\frac{1}{n^2} \right).
\end{equation}

Denote by $m_i$ the number of dispensable edges exposed during stage $i$. Then,
\begin{equation}
\label{eq:appendix:Bi}
    |B_{i+1}| \ge (d-1)|B_i| - 2m_{i+1} - m_i,
\end{equation}
where we used that every dispensable edge revealed at stage $i$ (resp. $i+1$) is incident to at most one vertex (resp. two vertices) in $B_i$.
Thus, after running the process for two full stages, by~\eqref{eq:appendix:1dis} and~\eqref{eq:appendix:Bi}, with probability $1 - O(n^{-2})$, we have $|B_1| \ge d-1$ and $|B_2| \ge (d-1)^2 - 2$. 
Moreover, denoting $B_1(x)$ and $B_2(x)$ for the same set with a vertex $x \in V$ as a starting point (i.e., we have $B_1 = B_1(v)$ and $B_2 = B_2(v)$), by taking a union bound over all $x \in V$, we get that, with probability $1 - O(n^{-1})$ we have $|B_1(x)| \ge d-1$ and $|B_2(x)| \ge (d-1)^2-2$ for all $x \in V$.
We assume in the sequel that this event holds.

We wish to show that, in the defined BFS process, the graph expands rather quickly. 
This is done by showing that at most one dispensable edge is revealed during the first $i=O(1)$ stages with suitably high probability and therefore, roughly speaking, the graph is very close to being a $(d-1)$-ary tree.
Indeed, even if dispensable edges occur after the very first few stages of the process, the graph has already expanded sufficiently so that these edges do not ``prevent'' further expansion in a significant way.

More formally, we bound from below the sizes of $B_i$.
By~\eqref{eq:appendix:Bi}, we can write
\begin{align}
\begin{split}
\label{eq:appendix:Lismall}
    |B_{i+1}| & \ge (d-1)^i |B_1| - \sum_{j=1}^{i} 2(d-1)^{i-j} m_{j+1} - \sum_{j=1}^{i} (d-1)^{i-j}m_j \\
    &\ge (d-1)^{i+1} - \sum_{j=1}^{i} 2(d-1)^{i-j} m_{j+1} - \sum_{j=1}^{i}(d-1)^{i-j}m_{j}.
\end{split}
\end{align}
As we encounter at most one dispensable edge during the first $i$ stages by assumption, we obtain that
$\sum_{j=1}^{i+1}m_{j} \le 1$ and, in particular,
\[\sum_{j=1}^{i} 2(d-1)^{i-j}m_{j+1} + \sum_{j=1}^{i}(d-1)^{i-j}m_j \le 2(d-1)^{i-1} + (d-1)^{i-2}. \]
Hence, for $i=O(1)$ (implying that $(d-1)^{i+1}=O(1)$), by~\eqref{eq:appendix:Lismall}, we obtain that 
\begin{align}\label{eq:B_Cbound}
    |B_{i+1}| \ge (d-1)^{i+1} - 2(d-1)^{i-1} - (d-1)^{i-2} \ge (d-9/4)(d-1)^i.
\end{align}

Now, consider the process after $k = o(n^{1/7})$ many steps in total (over all stages so far).
Similarly to~\eqref{eq:appendix:1dis}, note that the probability that more than two dispensable edges were revealed during the first $k = o(n^{1/7})$ steps in the process is at most
\begin{equation}
\label{eq:appendix:2dis}
    \binom{k}{3} \left(\frac{2k}{n-2k-1} \right)^3 = o\left(\frac{1}{n^{15/7}} \right).
\end{equation}
Fix an integer $i_0 = i_0(n)$ such that $n^{1/8} \le 2^{i_0+1} = o(n^{1/7})$.
Then, similarly to~\eqref{eq:appendix:Lismall}, for some large enough constant $C$ and for every $i=\omega(1)$ with $i = i(n) \le i_0$, we can write
\begin{align*}
    |B_{i+1}| &\ge (d-1)^{i+1-C}|B_C| - \sum_{j=C}^{i}2(d-1)^{i-j}m_{j+1} - \sum_{j=C}^{i} (d-1)^{i-j} m_j. 
\end{align*}
Recall that, with probability $1 - o(n^{-15/7})$, we have $\sum_{j=C}^{i_0+1} m_{j} \le 2$. Thus,
\begin{align*}
    \sum_{j=C}^{i} 2(d-1)^{i-j} m_{j+1} + \sum_{j=C}^{i} (d-1)^{i-j} m_j &\le 4(d-1)^{i-C} + 2(d-1)^{i-C-1}.
\end{align*}
Hence, by combining~\eqref{eq:B_Cbound} with the last two inequalities, we obtain that
\begin{align}
\begin{split}
\label{eq:appendix:Limed}
    |B_{i+1}| &\ge (d-1)^{i+1-C} \cdot (d-9/4)(d-1)^{C-1} - 6(d-1)^{i-C} \ge (d-7/3)(d-1)^{i}
\end{split}
\end{align}
holds with probability $1 - o(n^{-15/7})$ for all $i \le i_0$.

Let $t_0$ be the maximal integer such that
\begin{equation}
\label{eq:appendix:t0}
    (d-1)^{t_0} \le \sqrt{n} \log n.
\end{equation}
We now wish to show that $|B_{i+1}| \ge (d-5/2)(d-1)^{i}$ holds with probability $1 - o(n^{-17/8})$ for every $i\in [i_0, t_0]$.
Fix such an index $i$. Note that $S_i = \bigcup_{j=0}^{i} B_j$ and, moreover, that by the end of stage~$i+1$,~there were at most $d|S_i| = d\sum_{j=0}^{i} |B_j|$ steps in the process.
Recall that $|B_0|=1$ and trivially $|B_j| \le d(d-1)^{j-1}$ for all $j \ge 1$, so in particular $|S_i| \le 1+d+d\sum_{j=1}^{i}(d-1)^{j-1}$.
By~\eqref{eq:appendix:Pdis}, we obtain that the probability that an edge incident to a vertex in $S_i$ is dispensable is at most 
\begin{equation}
\label{eq:appendix:Si}
    \frac{d|S_i|}{n - d|S_i|} \le \frac{d+d^2+d^2\sum_{j=1}^{i}(d-1)^{j-1}}{n-d - d^2 - d^2\sum_{j=1}^{i}(d-1)^{j-1} } \le \frac{2d^2(d-1)^{i}}{n - 2d^2(d-1)^{i}} \le \frac{(d-1)^{t_0+5}}{n} \le \frac{(\log n)^2}{\sqrt{n}},
\end{equation}
for large enough $n$.
Let $X_i$ be the random variable counting the number of dispensable edges incident to $S_i$ in the entire process (that is, over the first $i+1$ stages).
By~\eqref{eq:appendix:Si}, we get that (the distribution of) $X_i$ is stochastically dominated by a binomial random variable $X \sim \text{Bin}(N,q)$ with $N = 2d^2(d-1)^i$ and $q = (\log n)^2 n^{-1/2}$.
Hence, the expected number of dispensable edges incident to vertices in $S_i$ for $i\le t_0$ satisfies
\[\mathbb E[X_i] \le \mathbb E[X] \le 2d^2(d-1)^i \cdot \frac{(\log n)^2}{\sqrt{n}} \le 4d^2 (\log n)^3, \]
assuming that $n$ is large enough.
By the Chernoff bound (\Cref{lem:Chernoff}), we obtain that
\begin{equation}
\label{eq:appendix:WiCrnff}
    \Prob\left(X_i \ge n^{1/16} \right) \le \Prob\left[X \ge n^{1/16} \right] \le e^{-n^{1/17}}.
\end{equation}
Thus, similarly to~\eqref{eq:appendix:Lismall}, we obtain that, with probability $1 - o(n^{-17/8})$, for all $i \in [i_0,t_0]$, we have
\begin{align}
\begin{split}
\label{eq:appendix:Lilarge}
    |B_{i+1}| &\ge (d-1)^{i+1-i_0} |B_{i_0}| - \sum_{j=i_0}^{i} 2(d-1)^{i-j} m_{j+1} -\sum_{j=i_0}^{i} (d-1)^{i-i_0} m_j \\
    &\ge (d-1)^{i+1-i_0}\cdot (d-7/3)(d-1)^{i_0-1} - 2(d-1)^{i-i_0}n^{1/16} - (d-1)^{i-i_0}n^{1/16} \\
    &\ge (d-5/2)(d-1)^{i}.
\end{split}
\end{align}

We now claim that, in this process, the graph expands well even after deleting up to $(\log n)^2$ many edges from it, uniformly at random.
We claim that at most two edges were deleted up to the end of stage $t_0$.
As we have $(d-1)^{t_0} \le \sqrt{n}\log n$, we get that this happens with probability at most
\begin{align*}
\binom{\sqrt{n}\log n}{3} \left(\frac{10(\log n)^2}{dn} \right)^3 = O\left(\frac{(\log n)^{9}}{n^{3/2}} \right) = o(n^{-4/3}).
\end{align*}
Taking a union bound over all vertices we get that with probability $1-o(n^{-1/3})$, for all vertices $x \in V$ we have that at most two edges were deleted in the graph obtained by running the process for at most $t_0$ many stages starting at $x$.
Similarly, we assume that this event holds for all vertices in the sequel.

Now note that, for every $v \in V$, there are $d \ge 3$ edges incident to it.
Hence, in the exploration process described above, after deleting edges, for all $v\in V$, there is at least one edge $vu$ incident to $v$ such that all the edges explored away from $u$ were not deleted.
In particular, a BFS exploration process after the deletion of edges stochastically dominates the original BFS process restricted to $B_1 = \{u\}$ after the first stage but without deletion of edges.
For $i \le t_0$, let $B'_{i}(x)$ be the set of vertices of distance precisely $i$ from $x$ after deleting edges as described.
Hence, we have $|B'_i(x)| \ge (d-5/2)(d-1)^{i-2}$.

To finish the proof of the statement, we consider two vertices $u,v \in V$.
Fix $s_1 = \lceil s/2\rceil$ and $s_2=\lfloor s/2\rfloor+1$ with $s$ defined in the statement of the theorem.
We now show that, with probability $1-o(n^{-17/8})$, either there is an edge between $B'_{i}(u)$ and $B'_{j}(v)$ for some $i\le s_1$ and $j\le s_2$ 
which survives after the deletion of edges by assumption, or there is an edge between $B'_{s_1}(u)$ and $B'_{s_2}(v)$ surviving the edge deletion.
To this end, we bound the following probability 
\begin{equation*}
    \Prob\left(\dist(u,v) > s ~|~ \dist(u,v) > s-1 \right) = \Prob\left(\dist(u,v) > s_1+s_2-1 ~|~ \dist(u,v) > s_1+s_2-2 \right).
\end{equation*}
Recall that, by the assumption of the statement, $s$ is the least integer satisfying $(d-1)^{s-1} \ge (4+\varepsilon)dn \log n$.
Hence, by~\eqref{eq:appendix:t0}, we have $s_1, s_2 \le t_0$.
By~\eqref{eq:appendix:Lismall}, \eqref{eq:appendix:Limed} and~\eqref{eq:appendix:Lilarge}, with probability $1-o(n^{-17/8})$, we obtain that $|B'_{s_1}(u)| \ge (d-5/2)(d-1)^{s_1-2}$ and $|B'_{s_2}(v)| \ge (d-5/2)(d-1)^{s_2-2}$.
In particular, with probability $1 - o(n^{-17/8})$, there are at least $(d-5/2)(d-1)^{s_1-2}$ edges with one endvertex in $B'_{s_1}(u)$ and another endvertex outside $B'_{s_1}(u)$: indeed, there are at most two dispensable edges in the graph $S_{s_1}(u)\supseteq B_{s_1}'(u)$. A similar statement holds for $B_{s_2}'(v)$.

At the same time, for every fixed half-edge emanating from $B_{s_1}'(u)$, the probability that at least one edge goes between $B'_{s_1}(u)$ and $B'_{s_2}(v)$ is at least $(d-5/2)(d-1)^{s_2-2}/(dn)$.
Denote the number of half-edges explored neither from $u$ after $s_1$ stages nor from $v$ after $s_2$ stages of the exploration process by $r$, and the number of half-edges emanating from $B'_{s_1}(u)$ and $B'_{s_2}(v)$ by $t_1$ and $t_2$, respectively.
Then, the probability of having at most 3 edges between $B'_{s_1}(u)$ and $B'_{s_2}(v)$ is bounded from above by
\[\sum_{i=0}^{3} \binom{t_1}{i} \bigg(\prod_{j=0}^{i-1} \frac{t_2-j}{r-2j-1}\bigg)\bigg(\prod_{j=0}^{t_1-i} \frac{r-(t_2+i)-2j-1}{r-2i-2j-1}\bigg)\le \sum_{i=0}^{3} \frac{2}{i!} \bigg(\frac{t_1t_2}{r}\bigg)^i \exp\bigg(-\sum_{j=0}^{t_1-i} \frac{t_2-i}{r-2i-2j-1}\bigg),\]
where the last inequality used that $1-x\le \e^{-x}$ for every $x\ge 0$.
Using that 
\[(4+\eps/4)dn\log n\le (d-5/2)^2(d-1)^{s-3}\le t_1t_2\le d^2(d-1)^{s-3}\le (4+\eps)2d^2n\log n.\]
and $r = dn - o(n)$, the exponential term above is of order $\exp(-(1+o(1)) t_1t_2/r) = o(n^{-4})$. As a consequence,
\begin{align*}
    \Prob\left(\dist(u,v) > s \right) \le \Prob\left(\dist(u,v) > s ~|~ \dist(u,v) > s-1 \right) = o(n^{-4}).
\end{align*}
By taking a union bound over all pairs of vertices, we obtain that
\begin{align*}
    \Prob\left(\diam(H) > s \right) &= \Prob\left(\exists u,v \in V : \dist_H(u,v) > s \right) \le \binom{n}{2} \Prob\left(\dist_H(u,v) > s \right) = o\left((\log n)^{-2} \right),
\end{align*}
finishing the proof.
\end{proof}

\end{document}